\documentclass[reqno,a4paper,10pt]{amsart}
\usepackage{amsfonts,amsthm,amsmath,latexsym,bm,graphics,indentfirst,color}
\usepackage[hypertex]{hyperref}
\usepackage{geometry}
\newtheorem{theorem}{Theorem}[section]
\newtheorem{lemma}[theorem]{Lemma}
\newtheorem{proposition}[theorem]{Proposition}
\newtheorem{corollary}[theorem]{Corollary}
\newtheorem{remark}[theorem]{Remark}

\newcommand{\e}{\varepsilon}
\newcommand{\D}{\displaystyle}
\newcommand{\F}[2]{\frac{#1}{#2}}

\newcommand{\BE}{\begin{equation}}
\newcommand{\BEN}{\begin{equation*}}
\newcommand{\EE}{\end{equation}}
\newcommand{\EEN}{\end{equation*}}
\newcommand{\BL}{\begin{lemma}}
\newcommand{\EL}{\end{lemma}}
\newcommand{\BT}{\begin{theorem}}
\newcommand{\ET}{\end{theorem}}
\newcommand{\BP}{\begin{proposition}}
\newcommand{\EP}{\end{proposition}}
\newcommand{\BC}{\begin{corollary}}
\newcommand{\EC}{\end{corollary}}
\newcommand{\BR}{\begin{remark}}
\newcommand{\ER}{\end{remark}}

\allowdisplaybreaks

\newcommand{\I}{\mathrm i}
\newcommand{\E}{\mathrm e}

\title[concentrations]{Concentrations for nonlinear Schr\"{o}dinger equation with magnetic potential and constant electric potential}

\author{Liping Wang}
\address{Liping Wang School of Mathematical Sciences, Shanghai Key Laboratory of Pure Mathematics and Mathematical Practice, East China Normal University, 200241, P.R. China}
\email{lpwang@math.ecnu.edu.cn}

\author{Chunyi Zhao}
\address{Chunyi Zhao School of Mathematical Sciences, Shanghai Key Laboratory of Pure Mathematics and Mathematical Practice, East China Normal University, 200241, P.R. China}
\email{cyzhao@math.ecnu.edu.cn}

\thanks{The two authors are partially supported by NSFC 11971169 and Natural Science Foundation of Shanghai 22ZR1421600.}



\begin{document}
\maketitle

\begin{abstract}
This paper studies the concentration phenomena to nonlinear Schr\"{o}dinger equations with magnetic potential and constant electric potential. The existing results show that the magnetic field has no effect on the location of point concentrations, as long as the electric potential is not constant. This paper finds out what the role of  the magnetic field plays  in the location of concentrations when the electric potential is constant.
%It is the first result in showing the role of the magnetic field with constant electric potential.
\end{abstract}

\section{Introduction}
%Assumption $\zeta = \zeta_0 + \e \xi$, $\partial_i A_i (\zeta_0) = \textcolor[rgb]{1.00,0.00,0.00}{\nabla\partial_i A_i (\zeta_0)}=0 $ \par
In this paper, we investigate the magnetic Schr\"odinger equation  in $\mathbb R^N$
\BE\label{tse}
\I\e \frac{\partial \psi}{\partial t}= (\I\e\nabla + \bm A(x))^2 \psi + Q(x) \psi - |\psi|^{p-1} \psi,
\EE
where $\e >0$ is a small parameter, $1<p<\frac{N+2}{N-2}$ if $N \ge 3$, $p>1$ if $N=2$ and $\I$ is the imaginary unit. The function $\psi$ is   complex-valued. Vector $\bm A=(A_1,A_2,\ldots,A_N): \mathbb R^N \to \mathbb R^N $ is the magnetic potential and assumed to be smooth, bounded and  $Q: \mathbb R^N \to \mathbb R$ represents the electric potential.  The magnetic Laplacian is defined by
\[
(\I\e\nabla + \bm A)^2 \psi :=-\e^2\Delta \psi +2\I\e \bm A \cdot \nabla \psi + |\bm A|^2 \psi +\I\e\psi \nabla\cdot \bm A.
\]
Terms involving $\bm A$ model the presence in some quantum model of the  magnetic field $\bm B$ given by
\BEN
\bm B=
\begin{cases}
\partial_1 A_2 - \partial_2 A_1, & \text{for }N=2,\\
\text{curl} \bm A, & \text{for }N=3, \\
(\partial_j A_k - \partial_k A_j)_{N\times N}, & \text{for }N >3.
\end{cases}
\EEN
For the discussion of this operator, one may refer to  \cite{L} and  \cite{M}.

The equation (\ref{tse}) arises in various physical contexts such as Bose-Einstein condensates and  nonlinear optics, see  \cite{BW}; or plasma physics where one can simulate the interaction effect among many particles by introducing some nonlinear term, see \cite{SS}. Concerning nonlinear Schr\"{o}dinger equation with the  magnetic field, the  pioneer  work is by Esteban-Lions \cite{EL} where they studied some minimization problems under suitable assumptions on
the magnetic field by concentrations and compactness arguments. For more results, one can refer to  \cite{AS, H, HE, K} and the references therein.

From now on we consider standing wave solutions to problem (\ref{tse}), namely $\psi(t,x)=\E^{\I\lambda \e^{-1}t}u(x)$ for some complex-valued function $u(x)$. Substituting this ansatz into problem (\ref{tse}),  $u(x)$ should satisfy
\BE\label{NLSE}
 (\I\e\nabla + \bm A)^2 u + V(x) u - |u|^{p-1} u =0 \qquad \text{in }\mathbb R^N,
\EE
where $V(x)=Q(x) + \lambda$. The potential $V(x)$ is usually assumed to be smooth and
\BEN
\inf_{\mathbb R^N} \ V(x) > 0.
\EEN
Problem (\ref{NLSE}) seems a very interesting problem since the Correspondence's Principle establishes that classical mechanics is, roughly speaking, contained in quantum mechanics.  The mathematical transition from quantum mechanics to classical mechanics can be formally described by letting the Planck constant  $\e \rightarrow 0$, and thus the existence of solutions for $\e$ small has physical interest. Standing waves for $\e$ small are usually referred as semi-classical bound states, see e.g. \cite{H}. In the linear case, Helffer et al. in \cite{HE, HS} have studied the asymptotic behavior of the eigenfunctions of the Schr\"{o}dinger operators with magnetic fields in the semiclassical limit. Note that in these results, one can not see the magnetic field in the definition of the {\it well}.

Concentration plays the important role in the problem (\ref{NLSE}). The pioneer work of concentration to problem (\ref{NLSE}) with $\bm A \equiv 0$ in one dimension is given by Floer and Weistein \cite{FW}. Really they proved that if the electric potential $V(x)$ has a non-degenerate critical point, then $u(x)$ concentrates near this critical point as $\e\rightarrow 0$.  After that, it appears many results concerning concentration phenomena about problem (\ref{NLSE}). Different approaches are used to cover different cases, see \cite{DFE, LI, OH, WZE}. Cingolani and Secchi \cite{CS1} proved that for bounded vector potentials,  concentration can happen at any non-degenerate critical point, not necessarily a minimum of $V$, as $\e \rightarrow 0$ using a perturbation approach given by Ambrosetti, Malchiodi and Secchi in \cite{AM}. Later they \cite{CS2} extended the result to the case of a possibly unbounded vector potential by a penalization procedure (see \cite{DF}). Semiclassical multi-peak solutions were found for bounded vector potentials in \cite{AB, CN, CT, GUI}.  Li-Peng-Wang \cite{LPW} constructed infinitely many  non-radial solutions if  the potentials $\bm A$ and $V$ are radially symmetric and suitably decay at infinity by using the number of bumps as parameter motivated by \cite{WY}, see also e.g.  \cite{WWY} and the references therein  for the application of such a very novel idea. For more recent results, we can refer to \cite{CC, CJS, T1, T2} and the references therein. Also readers can refer to \cite{AMN, DKW, WWY1, WZ} for high dimensional concentration. Till now, all these concentrations, especially point concentrations, are very dependent on   critical points of  $V(x)$ while the effect of the magnetic vector potential $\bm A$ is always ignored as higher order. In other words,  $\bm A$ has no contribution to decide the location of point concentrations when $V(x)$ has critical points. On the other hand,  for the case of the constant $V(x)$, there is no result on the concentration  as far as we know.    \par

Our aim in this paper is to exhibit how the magnetic potential $\bm A$ plays the role in the concentration. For this purpose, we study  the following nonlinear magnetic Schr\"odinger equation
\BE\label{main problem}
(\I \e \nabla + \bm A(x))^2 u + u - |u|^{p-1} u = 0 \qquad \mbox{in } \mathbb R^N.
\EE
Here $p>1$ and  $\bm A=(A_1,\ldots, A_N)$: $\mathbb R^N \to \mathbb R^N$ is assumed in $W^{1,\infty}(\mathbb R^N, \mathbb R^N)$ and smooth everywhere for simplicity.
%Obviously, for the constant vector $\bm A$,  the  problem (\ref{main problem}) has  a solution $w(\e^{-1}x)\E^{\I\e^{-1} \bm A\cdot x }$  where $w(y)=w(|y|)$ is the unique radial real-valued solution of
%\begin{equation}\label{rvp}
%\Delta w - w + w^p = 0 \quad \text{in }\mathbb R^N,  \qquad w(0) = \max_{\mathbb R^N} w >0, \qquad w(\pm\infty)=0.
%\end{equation}
To state our main result, we introduce the Frobenius norm of a matrix $\bm{\mathcal M} = (m_{ij})_{I\times J}$
\BEN
\|\bm{\mathcal M}\|_F = \left(\sum_{i=1}^I \sum_{j=1}^J |m_{ij}|^2\right)^\F{1}{2}.
\EEN
And denote by $w(y)=w(|y|)$  the unique radial real-valued solution of
\begin{equation}\label{rvp}
\Delta w - w + w^p = 0 \quad \text{in }\mathbb R^N,  \qquad w(0) = \max_{\mathbb R^N} w >0, \qquad w(\pm\infty)=0.
\end{equation}
Now our main results are the following.

\BT\label{t1}
Assume that the Forbenius norm $\|\bm B\|_F$ of the magnetic field  $\bm B$ admits $K$ local maximum (minimum) points  $\{P_m\}_{m=1}^K$ (which may be degenerate) and $K$ disjoint, closed and  bounded  regions $\{\Omega_m\}_{m=1}^K$ of $\mathbb R^N$ such that
\BEN
\|\bm B(P_m)\|_F=\max_{ \Omega_m} \|\bm B\|_F > \max_{\partial \Omega_m} \|\bm B\|_F. \qquad \left(\|\bm B(P_m)\|_F=\min_{ \Omega_m} \|\bm B\|_F < \min_{\partial \Omega_m} \|\bm B\|_F. \right)
\EEN
Then there exists an $\e_0>0$, such that for every $0<\e<\e_0$, the problem (\ref{main problem}) admits a  solution $u_\e$ with the form
\BEN
u_\e(x) =\sum_{m=1}^K \left[w\left(\F{|x-\zeta_{m}|}{\e}\right) + \e\Psi_m\left(\frac{x}\e\right)\right] \E^{\I\sigma_m+\I\e^{-1} \bm A(\zeta_{m})\cdot x} +  O(\e^2),
\EEN
for some $(\sigma_1, \ldots, \sigma_K) \in [0, 2\pi)^K,  \zeta_{m}\in \Omega_m$. The definition of $\Psi_m$ is given in (\ref{psi}).
%Moreover, $|u_\e|$ possesses exactly $K$ local maximum points $\zeta_{m,\e}$ in $\Omega_m$ respectively.
\ET

For general critical points, i.e. not local extremum points, we also have the following result.

\BT \label{mt}
Assume that $p>\F{3}{2}$ and $P_1, P_2, \ldots, P_K$, $K\ge 1$  are all non-degenerate critical points of $\|\bm B\|^2_F$.
Then there exists an $\e_0>0$ such that for any $ 0<\e <\e_0,$ the problem (\ref{main problem}) admits a solution $u_\e$ with the form
\BEN
u_\e(x) =\sum_{m=1}^K \left[w\left(\F{|x-\zeta_{m}|}{\e}\right) + \e\Psi_m\left(\frac{x}\e\right)\right] \E^{\I\sigma_m+\I\e^{-1} \bm A(\zeta_{m})\cdot x} +  O(\e^2)
\EEN
for some $(\sigma_1, \ldots, \sigma_K) \in [0, 2\pi)^K$ and   $\zeta_{m}=P_m + o(1)$. The definition of $\Psi_m$ is given in (\ref{psi}).
\ET

\begin{remark}
Theorem \ref{t1} and Theorem \ref{mt} clearly show the principle of  the magnetic vector $\bm A$, or precisely the magnetic field $\bm B$   driving the location of concentration if the electric potential is constant. This is a completely new phenomenon up to now.
\end{remark}

\begin{remark}
  Theorem \ref{t1} holds for any $p>1$. And Theorem \ref{mt} holds only for $p>\F{3}{2}$. The reason is that  the critical point of the function deduced by the energy functional is the local extremum point  in Theorem \ref{t1}, which may be gotten by direct comparison. While in Theorem \ref{mt}, $P_1, P_2, \ldots, P_K$ may  not be extremum points any more. Thus we have to study the derivatives of the energy functional. Naturally,  the corresponding estimates  should be higher accuracy  for which the better regularity of  the nonlinear term $|u|^{p-1}u$ is required.
\end{remark}

     From the point view of physics,  the magnetic field $\bm B$   is essential, not the particular choice of magnetic  potential $\bm A$. At the same time, there is a gauge invariance  for the magnetic Laplacian correspondingly, that is,  the magnetic field $\bm B$ is invariant under the transform of the potential $\bm A\to \bm A +\nabla f$. Also it is easy to see that the energy of the problem (\ref{main problem}) is unchanged under the gauge invariance, with which our result coincides (see Proposition \ref{p4.3}).

  When the electric potential $V(x)$ has some  non-degeneracy,   we just need to make the first approximation by the limit equation, which is enough to deduce the role of the function $V(x)$.  But in order to  see the role of magnetic vector potential $\bm A$ in our case, we need the  approximation up to order $\e$ of the problem (\ref{main problem}). Hence it is necessary to get a more accurate expression. On the other hand, higher accuracy allows us to deal with extremum points not necessary minimum points, which is different from the case of the non-constant $V(x)$, see e.g. \cite{GUI}.  Finally we note that the solutions given in Theorem \ref{t1} and Theorem \ref{mt} both show  simple concentrations. In the forthcoming paper, we will deal with multi-bump solutions for constant electric potential.

The paper is organized as follows. In section \ref{s2},  the ansatz and the estimate of the error are given.   The corresponding non-linear problem is solved  in section 3. In section 4,  the original problem is reduced to the finite dimension problem using variational reduction process  and  the expansion of the energy functional is shown. Finally, in section 5, Theorem \ref{t1} and Theorem \ref{mt} are proved.

\textbf{Notations.} \\
1. Constant $\beta=\min\{p-1,1\}$. \\
2. The real part of  $z \in \mathbb{C}$ will be denoted by $\operatorname{Re} z$. \\
3. The complex conjugate of  $z \in \mathbb{C}$ will be denoted by $\bar{z}$.\\
4. $C$ denotes a generic positive constant, which may be different from lines to lines. \\
5. Landau symbols $O(\e)$ is a generic function such that $|O(\e)| \le C\e$ and $o(\e)$ means that
 $\lim_{\e \rightarrow 0^+} \frac{|o(\e)|}{\e} =0$. \par

\section{Ansatz}\label{s2}

In this section we will present the approximation  of the problem and give the corresponding error estimate. \par

Recall in \cite{K} that the problem
\BEN
\Delta \tilde w  - \tilde w +|\tilde w|^{p-1} \tilde w =0, \qquad \tilde w \in H^1(\mathbb R^N, \mathbb C)
\EEN
possesses a unique ground state solution $\tilde w(y) = w(y) \E^{\I \sigma}$,  $\forall\ \sigma \in [0,2\pi]$ where $w(y)$ is the radial solution of problem (\ref{rvp}).    Thus  by the gauge invariance,
\BEN
\widetilde U(x) = w\left(\F{|x-\zeta|}{\e}\right) \E^{\I\sigma+\I\e^{-1} \bm A(\zeta)\cdot x}, \qquad \forall\ \sigma\in [0, 2\pi],
\EEN
is also the ground state solution to the constant magnetic potential problem
\BEN
(\I\e \nabla + \bm A(\zeta))^2 \widetilde U + \widetilde U - |\widetilde U|^{p-1} \widetilde U = 0, \qquad \widetilde U\in H^1(\mathbb R^N,\mathbb C).
\EEN

In the frame of large variable $y=x/\e$, the original problem (\ref{main problem}) is equivalent to
\BE\label{appro}
(\I \nabla + \bm A(\e y))^2 u + u - |u|^{p-1} u = 0 \qquad \mbox{in } \mathbb R^N.
\EE
Therefore, the function
%we still write for convenience
\BEN
U(y) = w(|y-\zeta'|)\E^{\I\sigma+\I \bm A(\zeta)\cdot y}, \qquad \zeta'=\zeta/\e.
\EEN
formally approximates the solution at least near $y=\zeta'$.
%Note that $U(y)$ is the solution of
%\BEN
%(\I \nabla + \bm A(\zeta))^2 U + U - |U|^2 U = 0 \qquad \text{in }\mathbb R^N.
%\EEN
Now it's time to give the first approximation
\begin{equation}
W(y)=\sum_{m=1}^K U_m(y), \qquad \quad U_m(y)= w(|y-\zeta_m'|)\E^{\I\sigma_m +\I \bm A(\zeta_m)\cdot y}, \qquad \sigma_m \in [0, 2\pi],
\end{equation}
where $\zeta_m \in \Omega_m, m=1,2,\ldots,K$. Denote
\begin{equation*}
\rho:=\min_{1\le i \neq j \le K}\text{dist}\left(\Omega_i,   \Omega_j \right)>0
\end{equation*}
and take a positive number $\delta$ such that  $ \delta\le\frac14\rho$.
%We aim to find an exact solution to (\ref{appro}) of the form $W( y)+\phi(y)$, where $\phi(y)$ is a small perturbation. So the function $\phi$ need to satisfy
%\BEN
%(\I \nabla + \bm A(\e y))^2 \phi + \phi - (p-1)|W|^{p-3}\text{Re}(\overline W \phi) W - |W|^2 \phi = - \widehat R(y) - N(\phi)
%\EEN
Let $\widehat{R}(y)$ be the error caused by the first approximation  $W$, which is
\BEN
\widehat R(y) = (\I \nabla + \bm A(\e y))^2 W + W -|W|^{p-1} W.
\EEN
 It is  checked that for $|y-\zeta'_1|\leq \frac{\delta}{\sqrt{\e}}, \quad |U_m|\leq |U_1| e^{-\frac{|\zeta'_m-\zeta'_1|}{3}}, m=2, 3,..., K$ and
\begin{align*}
\widehat{R}(y)&=\widetilde R(y) +O\left(w(|y-\zeta'_1|)\max_{2\le m\le K}e^{-\frac{|\zeta_m-\zeta_1|}{3\e}}  + w^{p-1}(|y-\zeta'_1|)\max_{2\le m\le K}e^{-\frac{|\zeta_m-\zeta_1|}{2\e}}\right)\\
&=\widetilde R(y) +O\left( w(|y-\zeta'_1|)+w^{ p-1}(|y-\zeta'_1|)\right)e^{-\frac{\delta}{\e}}
\end{align*}
where
\begin{align}
\widetilde R(y) &= (\I \nabla + \bm A(\e y))^2 U_1 + U_1 -|U_1|^{p-1} U_1 =  (\I \nabla + \bm A(\e y))^2 U_1 - (\I \nabla + \bm A(\zeta_1))^2 U_1  \nonumber \\
& = 2 \I (\bm A(\e y) -\bm A(\zeta_1))\cdot \nabla U_1 + \I \e (\nabla_x \cdot A) U_1  + \left(|\bm A(\e y)|^2 - |\bm A(\zeta_1)|^2\right) U_1  \nonumber \\
&= 2 \I (\bm A(\e y) -\bm A(\zeta_1))\cdot \left[\nabla w (y-\zeta'_1)+ \I \bm A(\zeta_1) w(y-\zeta'_1)\right]\E^{\I\sigma_1 +\I \bm A(\zeta_1)\cdot y}  + \I \e (\nabla_x \cdot \bm A) U_1 \nonumber\\
 &\quad + \left(|\bm A(\e y)|^2 - |\bm A(\zeta_1)|^2\right) w(y-\zeta'_1)\E^{\I\sigma_1 +\I \bm A(\zeta_1)\cdot y}  \nonumber \\
&= \left[-2 \bm A(\zeta_1)\cdot(\bm A(\e y) -\bm A(\zeta_1)) + \left(|\bm A(\e y)|^2 - |\bm A(\zeta_1)|^2\right)\right ] U_1 \nonumber\\
&\quad + \I [2 (\bm A(\e y) -\bm A(\zeta_1))\cdot\nabla w (y-\zeta'_1)+\e (\nabla_x \cdot A) w(y-\zeta'_1) ]\E^{\I\sigma_1+\I \bm A(\zeta_1)\cdot y}  \nonumber \\
&= \left(|\bm A(\e y)- \bm A(\zeta_1)|\right)^2 w(y-\zeta'_1)\E^{\I\sigma_1 +\I \bm A(\zeta_1)\cdot y}\nonumber\\
&\quad + \I \left[2 (\bm A(\e y) -\bm A(\zeta_1))\cdot\nabla w (y-\zeta'_1)+\e (\nabla_x \cdot \bm A) w(y-\zeta'_1) \right]\E^{\I\sigma_1+\I \bm A(\zeta_1)\cdot y}. \label{1}
\end{align}
As usual one writes $\zeta'_m=(\zeta'_{m,1}, \zeta'_{m,2}, \ldots, \zeta'_{m,N})$. Direct calculation shows that
\begin{align*}
& |\bm A(\e y)- \bm A(\zeta_1)|^2 = \sum_{i=1}^N (A_i(\e y) - A_i(\zeta_1))^2 \\
=&\ \sum_{i=1}^N \left(\e \nabla A_i(\zeta_1)\cdot (y-\zeta'_1) + \F{\e^2}{2} (y-\zeta'_1)^\perp\cdot \nabla^2 A_i(\zeta_1)\cdot(y-\zeta'_1) +O(\e^3 |y-\zeta'_1|^3)\right)^2 \\
%=&\ \e^2 \sum_{i,j,k=1}^N  \partial_j A_i(\zeta_1) \partial_k A_i(\zeta_1) (y_j-\zeta'_{1,j}) (y_k-\zeta'_{1,k})  \\
%&\qquad + \e^3 \sum_{i=1}^N \nabla A_i(\zeta_1)\cdot (y-\zeta'_1)(y-\zeta'_1)^\perp\cdot \nabla^2 A_i(\zeta_1)\cdot(y-\zeta'_1) + O(\e^4 |y-\zeta'_1|^4) \\
=&\ \e^2 \sum_{i,j,k=1}^N  \partial_j A_i(\zeta_1) \partial_k A_i(\zeta_1) (y_j-\zeta'_{1,j}) (y_k-\zeta'_{1,k})  \\
& + \e^3 \sum_{i,j,k,\ell=1}^N  \partial_j A_i(\zeta_1) \partial_{k\ell} A_i(\zeta_1)(y_j-\zeta'_{1,j}) (y_k-\zeta'_{1,k})  (y_\ell-\zeta'_{1,\ell})+ O(\e^4 |y-\zeta'_1|^4).
\end{align*}
Similarly, it is verified that
\begin{align*}
& (\bm A(\e y)- \bm A(\zeta_1)) \cdot\nabla w(|y-\zeta_1'|) =\sum_{i=1}^N( A_i(\e y)- A_i(\zeta_1)) \partial_i w(|y-\zeta_1'|) \\
=&\ \e \sum_{i,j=1}^N   \partial_j A_i(\zeta_1) (y_j-\zeta'_{1,j}) \partial_i w(|y-\zeta_1'|)  + \F{\e^2}{2} \sum_{i,j,k=1}^N   \partial_{jk} A_i(\zeta_1)(y_j-\zeta'_{1,j}) (y_k-\zeta'_{1,k})\partial_i w(|y-\zeta_1'|)   \\
& + \F{\e^3}{6}\sum_{i,j,k,\ell=1}^N \partial_{jk\ell} A_i(\zeta_1)(y_j-\zeta'_{1,j}) (y_k-\zeta'_{1,k}) (y_\ell-\zeta'_{1,\ell})\partial_i w(|y-\zeta_1'|)   + O\left(\e^4 |y-\zeta_1'|^4 |\nabla w(|y-\zeta_1'|) |\right).
%=&\ \e \sum_{i=1}^N \partial_i A_i(\zeta) (y_i-\zeta'_i) \partial_i w +\e \sum_{i=1}^N \sum_{j\neq i}  \partial_j A_i(\zeta) (y_j-\zeta'_j) \partial_i w + \F{\e^2}{2} \sum_{i,j,k}  \partial_{jk} A_i(\zeta)(y_j-\zeta'_j) (y_k-\zeta'_k)\partial_i w \\
%&\quad + \F{\e^3}{6}\sum_{i,j,k,\ell} \partial_{jk\ell} A_i(\zeta)(y_j-\zeta'_j) (y_k-\zeta'_k) (y_\ell-\zeta_\ell')\partial_i w  + O(\e^4 |y-\zeta'|^4 |\nabla w|) .
%=&\ \e \sum_{i=1}^N \sum_{j\neq i}  \partial_j A_i(\zeta) (y_j-\zeta'_j) \partial_i w + \e^{2} \sum_{i} \nabla(\partial_i A_i(\zeta_0))\xi (y_i-\zeta'_i) \partial_i w  \\
%&\quad + \e^2 \sum_{i} \sum_{j,k} \F{1}{2} \partial_{jk} A_i(\zeta)(y_j-\zeta'_j) (y_k-\zeta'_k)\partial_i w + O(\e^3 |y-\zeta'|^3 |\nabla w|)+ O(\e^{3} |\xi|^2 |y-\zeta'| |\nabla w|),
\end{align*}
As for the term $\nabla_x\cdot \bm A(\e y)w(|y-\zeta_1'|), $ we expand it similarly to
\begin{align*}
& \nabla_x\cdot \bm A(\e y)w(|y-\zeta_1'|)  \\
=& \nabla_x\cdot \bm A(\zeta_1)w(|y-\zeta_1'|)  + \e \nabla(\nabla_x\cdot \bm A)(\zeta_1)\cdot (y-\zeta_1')w(|y-\zeta_1'|)\\
& +\frac12\e^2(y-\zeta'_1)^\perp\cdot\nabla^2(\nabla_x\cdot \bm A)(\zeta_1)\cdot(y-\zeta'_1) + O\left(\e^3|y-\zeta'_1|^3 w(|y-\zeta_1'|) \right)\\
=& \nabla_x\cdot \bm A(\zeta_1)w(|y-\zeta_1'|)  +  \e  \sum_{i,j=1}^N  \partial_{ij}A_i(\zeta_1)(y_j-\zeta'_{1,j})w(|y-\zeta_1'|) \\
&  + \F{\e^2}{2} \sum_{i,j,k=1}^N \partial_{jki} A_i(\zeta_1) (y_j-\zeta'_{1,j}) (y_k-\zeta'_{1,k}) w(|y-\zeta_1'|)  + O(\e^3|y-\zeta_1'|^3 w(|y-\zeta_1'|) ).
%&= \e  \sum_{i}  \nabla \partial_i A_i(\zeta_0) \xi w + \e  \sum_{i,j}  \partial_{ij}A_i(\zeta)(y_j-\zeta_j')w + O(\e^{2}|\xi|^2 w ) + O(\e^2|y-\zeta'|^2 w).
\end{align*}
Therefore, in the region $|y-\zeta_1'| < \frac{\delta}{\sqrt{\e}}$, we get that
\begin{align*}
 &\ \widetilde R(y)\E^{-\I\sigma_1-\I \bm A(\zeta_1)\cdot y}\\
 =&\  \e \I \left( 2 \sum_{i,j=1}^N   \partial_j A_i(\zeta_1) (y_j-\zeta'_{1,j}) \partial_i w (|y-\zeta'_1|) + \nabla_x \cdot \bm A(\zeta_1) w(|y-\zeta'_1|) \right) \\
  & +\e^2 \sum_{i,j,k=1}^N \partial_j A_i(\zeta_1) \partial_k A_i(\zeta_1) (y_j-\zeta'_{1,j}) (y_k-\zeta'_{1,k}) w(|y-\zeta'_1|) \\
   & + \e^2 \I \sum_{i,j,k=1}^N   \partial_{jk} A_i(\zeta_1)(y_j-\zeta'_{1,j}) (y_k-\zeta'_{1,k})\partial_i w(|y-\zeta'_1|) + \e^2 \I\sum_{i,j=1}^N  \partial_{ij}A_i(\zeta_1)(y_j-\zeta'_{1,j})w(|y-\zeta'_1|) \\
& +  \e^3  \sum_{i,j,k,\ell=1}^N  \partial_j A_i(\zeta_1) \partial_{k\ell} A_i(\zeta_1)(y_j-\zeta'_{1,j}) (y_k-\zeta'_{1,k})  (y_\ell-\zeta'_{1,\ell})w(|y-\zeta'_1|) \\
& + \F{\e^3}{3} \I\sum_{i,j,k,\ell=1}^N \partial_{jk\ell} A_i(\zeta_1)(y_j-\zeta'_{1,j}) (y_k-\zeta'_{1,k}) (y_\ell-\zeta'_{1,\ell})\partial_i w (|y-\zeta'_1|)\\
&  + \F{\e^3}{2} \I \sum_{i,j,k=1}^N \partial_{jki} A_i(\zeta_1) (y_j-\zeta'_{1,j}) (y_k-\zeta'_{1,k}) w(|y-\zeta'_1|) \\
&  +   O\left(\e^4 |y-\zeta'_1|^4 | w(|y-\zeta'_1|)|\right) + \I O\left(\e^4 |y-\zeta'_1|^4 |\nabla w(|y-\zeta'_1|)| + \e^4 |y-\zeta'_1|^3 |w(|y-\zeta'_1|)|\right).
\end{align*}

Note that the imaginary part of $\widetilde R(y)\E^{-\I\sigma_1-\I \bm A(\zeta_1)\cdot y}$ is $ O(\e)$. It is of less accuracy for later application.   Complying with the guideline that the better approximation, the more possibility to get a  solution, we should improve the accuracy of the approximation. To this purpose,  we find that the real-valued function $\Psi_{1,ij}(y) =  \F{1}{2}(y_i-\zeta'_{1,i})(y_j-\zeta'_{1,j}) w(|y-\zeta'_1|)$ $(i\neq j)$ is the solution to
\BEN
-\Delta \Psi_{1,ij} + \Psi_{1,ij} -  w^{p-1}(|y-\zeta'_1|)\Psi_{1,ij} =-2 (y_j-\zeta'_{1,j}) \partial_i w (|y-\zeta'_1|) =-2 (y_i-\zeta'_{1,i}) \partial_j w (|y-\zeta'_1|).
\EEN
Besides, $\Psi_{1,ii}(y) =  \F{1}{2} (y_i-\zeta'_{1,i})^2 w(|y-\zeta'_1|) $ satisfies the equation
\BEN
-\Delta \Psi_{1,ii} + \Psi_{1,ii} -  w^{p-1}(|y-\zeta'_1|) \Psi_{1,ii} = -2 (y_i-\zeta'_{1,i}) \partial_i w (|y-\zeta'_1|) - w .
\EEN
%Then $\varphi(y)=\I \partial_j A_i(\zeta)\Psi_{ij} \E^{\I \bm A(\zeta)\cdot y} $ fulfills
%\BEN
%(\I \nabla + \bm A(\zeta))^2 \varphi + \varphi  - 2\operatorname{Re}(\overline U \varphi) U - |U|^2 \varphi = 2 \I \partial_j A_i(\zeta)(y_j-\zeta'_j) \partial_i w \E^{\I \bm A(\zeta)\cdot y}.
%\EEN
Then obviously the function
\begin{eqnarray*}%\label{Psi}
\Psi_1(y) &=&  \I \left(\sum_{i=1}^N \sum_{j\neq i} \partial_j A_i(\zeta_1)\Psi_{1,ij}(y) + \sum_{i=1}^N \partial_i A_i(\zeta_1)  \Psi_{1,ii}(y) \right)\E^{\I\sigma_1+\I \bm A(\zeta_1)\cdot y} \\
&=& \I \left(\sum_{i,j=1}^N \partial_j A_i(\zeta_1)\Psi_{1,ij}(y) \right)\E^{\I\sigma_1+\I \bm A(\zeta_1)\cdot y}: =\I \psi_1(y)\E^{\I\sigma_1+\I \bm A(\zeta_1)\cdot y}
\end{eqnarray*}
satisfies
\begin{eqnarray*}
&&(\I \nabla + \bm A(\zeta_1))^2 \Psi_1 + \Psi_1  - (p-1)|U_1|^{p-3}\operatorname{Re}(\overline U_1 \Psi_1) U_1 - |U_1|^{p-1} \Psi_1 \\
=&& \I\left[ - 2   \sum_{i,j=1}^N\partial_j A_i(\zeta_1) (y_j-\zeta'_{1,j}) \partial_i w(|y-\zeta'_1|)  -  \nabla_x \cdot \bm A(\zeta_1)w(|y-\zeta'_1|)\right]  \E^{\I\sigma_1+\I \bm A(\zeta_1)\cdot y},
\end{eqnarray*}
since $\text{Re}(\overline U_1 \Psi_1) = 0$.
%Besides, we also set $\Phi_i(y) = - \F{1}{2} (y-\zeta_i')^2 w(|y-\zeta'|) $ and
%\BEN
%\Phi(y) =  -\I \sum_i \partial_i A_i (\zeta) \Phi_i(y) \E^{\bm A(\zeta)\cdot y}.
%\EEN
%Then
%\BEN
%(\I \nabla + \bm A(\zeta))^2 \Phi + \Phi  - 2\operatorname{Re}(\overline U \Phi) U - |U|^2 \Phi = - 2 \I  \sum_{i=1}^N  \partial_i A_i(\zeta) (y_i-\zeta'_i) \partial_i w \E^{\I \bm A(\zeta)\cdot y} - \I \sum_i \partial_i A_i(\zeta) w(|y-\zeta'|).
%\EEN
%It is checked that the function $ w(r) + r w'(r)$ satisfies
%\BEN
%\textcolor[rgb]{1.00,0.00,0.00}{\Delta \Psi_2 -\Psi_2 + 3 w^2 \Psi_2 = 2w}
%\EEN
%Thus we set
%\BEN
%\Psi_2(y) =  \F{1}{2}(\nabla_x \cdot \bm A)(\zeta) \left[w(y-\zeta')+|y-\zeta'|w'(y-\zeta')\right] \E^{\I \bm A(\zeta)\cdot y},
%\EEN
%which means
%\BEN
%(i\nabla + \bm A(\zeta))^2 \Psi_2 + \Psi_2 - 2\operatorname{Re}(\overline U \Psi_2) U - |U|^2 \Psi_2 = -(\nabla_x \cdot \bm A(\zeta)) w(y-\zeta') \E^{\I \bm A(\zeta)\cdot y}
%\EEN
%Let $\Psi_1$ be the solution of
%\BEN
%(i\nabla + \bm A(\zeta))^2 \Psi_1 + \Psi_1 - 2\operatorname{Re}(\overline U \Psi_1) U - |U|^2 \Psi_1 = -(y-\zeta')\nabla \bm A(\zeta) \nabla w(y-\zeta') \E^{\I \bm A(\zeta)\cdot y}
%\EEN
%Denote $\Psi(y) =  (\Psi_1 + \Psi_2)$.
Moreover, the same computation may also be carried out in the region $|y-\zeta_m'|<\F{\delta}{\sqrt{\e}}$ for $m=2,\cdots,K$.
The ultimate approximation  is then selected as
\begin{equation}\label{ba}
\mathcal W(y) = W(y) + \e  \Psi(y),
\end{equation}
where $\Psi(y)=\sum_{m=1}^K \Psi_m(y)$. Here
\begin{equation}\label{psi}
 \Psi_m(y)=\I \left(\sum_{i,j=1}^N \partial_j A_i(\zeta_m)\Psi_{m,ij}(y) \right)\E^{\I\sigma_m+\I \bm A(\zeta_m)\cdot y}: =\I \psi_m(y)\E^{\I\sigma_m+\I \bm A(\zeta_m)\cdot y}
\end{equation}
and
\[
\Psi_{m,ij}(y)=\F{1}{2}(y_i-\zeta'_{m,i})(y_j-\zeta'_{m,j}) w(|y-\zeta'_m|).
\]

 Now the approximation  $\mathcal W$ is  good enough to help us find a solution. Precisely our aim is  to find a  solution with the form $\mathcal W(y) + \phi(y)$ of problem (\ref{main problem}), where $\phi$ is a small perturbation and  satisfies the equation
\BE\label{phi}
L\phi := (\I \nabla + \bm A(\e y))^2 \phi + \phi - (p-1) |\mathcal W|^{p-3}\operatorname{Re}(\overline{\mathcal W}\phi)\mathcal W - |\mathcal W|^{p-1}\phi = - R(y) + N(\phi).
\EE
Here $R(y)$ denotes  the error caused by $\mathcal W$, which is
\begin{gather*}
R(y) = (\I \nabla  + \bm A(\e y))^2 \mathcal W + \mathcal W - |\mathcal W|^{p-1}\mathcal W,
\intertext{and the nonlinear term}
N(\phi) = |\mathcal W+\phi|^{p-1}(\mathcal W +\phi)  - |\mathcal W|^{p-1}\mathcal W -(p-1)|\mathcal W|^{p-3}\operatorname{Re}(\overline{\mathcal W}\phi)\mathcal W - |\mathcal W|^{p-1}\phi.
\end{gather*}
Obviously, with the notation $\beta=\min\{ p-1, \ 1\}$,
\BEN
|N(\phi)| \leq C|\phi|^{1+\beta}.
\EEN
%
%In what follows, let us denote $R(y)$ be  the error caused by $\mathcal W$, that is
%\BEN
%R(y) = (\I \nabla  + \bm A(\e y))^2 \mathcal W + \mathcal W - |\mathcal W|^{p-1}\mathcal W.
%\EEN

To solve problem (\ref{phi}), it is important to estimate the error $R(y)$.

\BP\label{p2.1}
We have that
\begin{gather*}
\|R(y)\|_{L^2} \leq C\e^{2},
\intertext{and}
\|\partial_{\zeta_{m,k}'}R(y)\|_{L^2} \leq C\e^{2}, \qquad \|\partial_{\sigma_m}R(y)\|_{L^2} \leq C\e^{2}, \qquad \forall\ m=1,\ldots, K, \ k=1,\ldots, N.
\end{gather*}
\EP

\begin{proof}
First we consider the domain $|y-\zeta'_m|\le \frac\delta{\sqrt{\e}},$ without loss of generality, say $m=1$. Then
\begin{equation}
\begin{array}{lll}
R(y) &=&
\widehat{R}(y)  + \e  \left[(\I \nabla  + \bm A(\e y))^2 \Psi +   \Psi - (p-1)|W|^{p-3}\operatorname{Re}(\overline W \Psi) W -|W|^{p-1}\Psi \right]\nonumber \\
\nonumber\\
&& - \left[|W+\e \Psi|^{p-1}(W+\e\Psi) - |W|^{p-1}W - \e(p-1)|W|^{p-3}\text{Re}(\overline W \Psi) W - \e |W|^{p-1}\Psi  \right] \nonumber \\
\nonumber\\
&=& \widetilde R(y) -  \e\I \left(2   \sum_{i,j=1}^N  \partial_j A_i(\zeta_1) (y_j-\zeta'_{1,j}) \partial_i w(y-\zeta'_1) + (\nabla_x\cdot \bm A)(\zeta_1) w(|y-\zeta'_1|) \right)\E^{\I\sigma_1+\I \bm A(\zeta_1)\cdot y} \nonumber\\
\nonumber\\
&& + \e\left[(\I \nabla  + \bm A(\e y))^2 \Psi_1  -   (\I \nabla  + \bm A(\zeta_1))^2 \Psi_1 \right] + \e  \left[(\I \nabla  + \bm A(\e y))^2\sum_{m=2}^K \Psi_m + \sum_{m=2}^K  \Psi_m\right] \nonumber \\
\nonumber\\
&&  - \F{p-1}{2}\e^2 \psi_1^2(y) w^{p-2}(|y-\zeta'_1|)   \E^{\I\sigma_1+\I \bm A(\zeta_1)\cdot y} -\I \F{p-1}{2}\e^3 \psi^3_1(y) w^{p-3}(|y-\zeta'_1|)  \E^{\I\sigma_1+\I \bm A(\zeta_1)\cdot y} \nonumber \\
\nonumber\\
&&+ O \left(\e^{4} |y-\zeta_1'|^8|w^{p}(|y-\zeta'_1|) + \E^{-\delta/\e} w(|y-\zeta_1'|) +\E^{-\frac{\delta}{\e}}w^{p-1}(|y-\zeta'_1|) \right).
\end{array}
\end{equation}
%where the following facts are used that
%\[
%\left||y-\zeta'_m|^2U_mU_1\right| + |\Psi_m\Psi_1| \le Ce^{-\frac\delta\e}w(y-\zeta'_1), \qquad  2\le m \le K.
%\]
Obviously we have that
\begin{eqnarray*}
&&\e^2 \psi_1(y)^2 w^{p-2}(y-\zeta'_1) = \e^2 w^{p-2}(|y-\zeta'_1|)\sum_{i,j=1}^N  \partial_j A_i(\zeta_1)\Psi_{1,ij}(y) \sum_{k,\ell=1}^N  \partial_\ell A_k(\zeta_1)\Psi_{1,k\ell}(y) \\
=&& \F{\e^2}{4} w^{p}(|y-\zeta'_1|)\sum_{i,j,k,\ell=1}^N \partial_j A_i(\zeta_1)\partial_\ell A_k(\zeta_1)(y_i-\zeta'_{1,i})(y_j-\zeta'_{1,j}) (y_k-\zeta'_{1,k})(y_\ell-\zeta'_{1,\ell})
\end{eqnarray*}
and
\[
\e  \left|\left[(\I \nabla  + \bm A(\e y))^2\sum_{m=2}^K \Psi_m + \sum_{m=2}^K  \Psi_m\right]\right|=O\left(\sum_{m=2}^K w(y-\zeta'_m)   \right)=O\left( \E^{-\frac\delta\e}w(y-\zeta'_1)   \right).
\]
Moreover, it is checked as in (\ref{1}) that
\begin{align}
 &\ (\I \nabla  + \bm A(\e y))^2 \Psi_1  -   (\I \nabla  + \bm A(\zeta_1))^2 \Psi_1 \nonumber \\
 =&\ 2\I (\bm A(\e y)- \bm A(\zeta_1))\cdot\nabla \Psi_1 + \I \e (\nabla_x\cdot\bm A(\e y))\Psi_1 + (|\bm A(\e y)|^2 - |\bm A(\zeta_1)|^2)\Psi_1  \nonumber \\
 =&\ -\left[2(\bm A(\e y)-\bm A(\zeta_1))\cdot\nabla \psi_1 + \e(\nabla_x \cdot \bm A(\e y)) \psi_1 \right] e^{\I\sigma_1+\I \bm A(\zeta_1)\cdot y} +  \I |\bm A(\e y)-\bm A(\zeta_1)|^2  \psi_1 e^{\I\sigma_1 +\I \bm A(\zeta_1)\cdot y}. \label{2}
\end{align}
Let us estimate them one by one in (\ref{2}) for $y\neq \zeta'_1$.
Note that
\begin{align*}
  &\ (\bm A(\e y)-\bm A(\zeta_1))\cdot\nabla \psi_1 = \sum_{k=1}^N (A_k(\e y)-A_k(\zeta_1)) \partial_k \psi_1  \\
   =&\   \e \sum_{k=1}^N \nabla  A_k(\zeta_1)\cdot(y-\zeta'_1)\partial_k \psi_1
   + \F{\e^2}{2} \sum_{k=1}^N (y-\zeta'_1)^\perp \cdot \nabla^2 A_k(\zeta_1)\cdot(y-\zeta'_1) \partial_k \psi_1 + O(\e^3 |y-\zeta'_1|^3 |\nabla \psi_1|)  \nonumber \\
  =&\ \e \sum_{k,\ell=1}^N \partial_\ell A_k(\zeta_1) (y_\ell-\zeta'_{1,\ell}) \left[ \sum_{i,j=1}^N \partial_j A_i(\zeta_1)\partial_k \Psi_{1,ij} \right] \\
  &\ + \F{\e^2}{2} \sum_{k,\ell, s=1}^N \partial_{\ell s}A_k(\zeta_1) (y_\ell-\zeta'_{1,\ell})(y_s-\zeta'_{1,s})\left[ \sum_{i,j=1}^N\partial_j A_i(\zeta_1)\partial_k \Psi_{1,ij} \right]   + O\left(\e^3 |y-\zeta'_1|^5 w(|y-\zeta'_1|)\right). \\
\end{align*}
It is easy to see  that ($\delta_{ik}$ is the Kronecker symbol) for $y\neq \zeta'_1$,
\begin{align*}
\partial_k \Psi_{1, ij}(y) &=  \F{1}{2}\left [\left(\delta_{ik} (y_j-\zeta'_{1,j}) +\delta_{jk} (y_i-\zeta'_{1,i})\right)w(|y-\zeta'_1|)+ (y_i-\zeta'_{1,i})(y_j-\zeta'_{1,j})(y_k-\zeta'_{1,k})\F{w'(|y-\zeta'_1|)}{|y-\zeta'_1|} \right].
%\partial_k \widetilde \Psi_i (y) &=- \F{1}{2}\left [ 2\delta_{ik}(y_i-\zeta_i') w(|y-\zeta'|) + (y_i -\zeta_i')^2(y_k -\zeta_k') \F{w'(|y-\zeta'|)}{|y-\zeta'|} \right].
\end{align*}
Thus we conclude that
\begin{align}
&\ (\bm A(\e y)-\bm A(\zeta_1))\cdot\nabla \psi_1 \nonumber \\
=&\ \F{\e}{2}\sum_{j,k,\ell=1}^N \partial_\ell A_k(\zeta_1) \left[ \partial_j A_k(\zeta_1)+ \partial_k A_j(\zeta_1) \right] (y_j-\zeta'_{1,j})(y_\ell-\zeta'_{1,\ell})w(|y-\zeta'_1|) \nonumber  \\
&\ + \F{\e}{2} \sum_{i,j,k,\ell=1}^N \partial_\ell A_k(\zeta_1)\partial_j A_i(\zeta_1)(y_i-\zeta'_{1,i})(y_j-\zeta'_{1,j})(y_k-\zeta'_{1,k})(y_\ell-\zeta'_{1,\ell})\F{w'(|y-\zeta'_1|)}{|y-\zeta'_1|} \nonumber \\
&\ + \F{\e^2}{4} \sum_{j, k,\ell,s=1}^N \partial_{\ell s}A_k(\zeta_1) \left[ \partial_j A_k(\zeta_1) + \partial_k A_j(\zeta_1)\right]  (y_j-\zeta'_{1,j})(y_\ell-\zeta'_{1,\ell})(y_s-\zeta'_{1,s}) w(y-\zeta'_1)  \nonumber \\
&\ + \F{\e^2}{4} \sum_{i,j,k,\ell,s=1}^N  \partial_{\ell s}A_k(\zeta_1) \partial_j A_i(\zeta_1) (y_i-\zeta'_{1,i}) (y_j-\zeta'_{1,j})(y_k-\zeta'_{1,k}) (y_\ell-\zeta'_{1,\ell})(y_s-\zeta'_{1,s})\F{w'}{|y-\zeta'_1|} \nonumber \\
&\ + O\left(\e^3 |y-\zeta'_1|^5 w(|y-\zeta'_1|)\right). \label{3}
\end{align}
Also it can be obtained  that
\begin{eqnarray*}
  (\nabla_x \cdot \bm A(\e y)) \psi_1 &=&   \F{1}{2}(\nabla_x \cdot \bm A)(\zeta_1)\sum_{i,j=1}^N  \partial_j A_i(\zeta_1)(y_i-\zeta'_{1,i})(y_j-\zeta'_{1,j})w(|y-\zeta'_1|) \\
  && + \F{\e}{2} \sum_{i,j,k,\ell=1}^N \partial_j A_i(\zeta_1) \partial_{\ell k} A_\ell(\zeta_1) (y_i-\zeta'_{1,i}) (y_j-\zeta'_{1,j}) (y_k-\zeta'_{1,k})  w (y-\zeta'_1)\\
  &&+ O\left(\e^2 |y-\zeta'_1|^4 w(y-\zeta'_1)\right),
\end{eqnarray*}
and
\BEN
|\bm A(\e y)-\bm A(\zeta_1)|^2  \psi_1 = O\left(\e^2 |y-\zeta'_1|^4 w(y-\zeta'_1)\right).
\EEN
Now (\ref{2}) can be expanded as
\begin{align*}
   & (\I \nabla  + \bm A(\e y))^2 \Psi_1  -   (\I \nabla  + \bm A(\zeta_1))^2 \Psi_1 \nonumber \\
   =& -\e \sum_{j,k,\ell=1}^N \partial_\ell A_k(\zeta_1) \left[ \partial_j A_k(\zeta_1)+ \partial_k A_j(\zeta_1) \right] (y_j-\zeta'_{1,j})(y_\ell-\zeta'_{1,\ell})w(|y-\zeta'_1|) e^{\I\sigma_1+\I \bm A(\zeta_1)\cdot y} \nonumber \\
  & - \e \sum_{i,j,k,\ell=1}^N \partial_\ell A_k(\zeta_1)\partial_j A_i(\zeta_1)(y_i-\zeta'_{1,i})(y_j-\zeta'_{1,j})(y_k-\zeta'_{1,k})(y_\ell-\zeta'_{1,\ell})\F{w'(|y-\zeta'_1|)}{|y-\zeta'_1|} \E^{\I\sigma_1+\I \bm A(\zeta_1)\cdot y} \\
  & -\F{\e}{2}(\nabla_x \cdot \bm A(\zeta_1))\sum_{i,j=1}^N  \partial_j A_i(\zeta_1)(y_i-\zeta'_{1,i})(y_j-\zeta'_{1,j})w(|y-\zeta'_1|)e^{\I\sigma_1+\I \bm A(\zeta_1)\cdot y} \\
  & - \F{\e^2}{2} \sum_{j, k,\ell,s=1}^N \partial_{\ell s}A_k(\zeta_1) \left[ \partial_j A_k(\zeta_1) + \partial_k A_j(\zeta_1)\right]  (y_j-\zeta'_{1,j})(y_\ell-\zeta'_{1,\ell})(y_s-\zeta'_{1,s}) w(y-\zeta'_1) \E^{\I\sigma_1+\I \bm A(\zeta_1)\cdot y} \\
  & - \F{\e^2}{2} \sum_{i,j,k,\ell,s=1}^N  \partial_{\ell s}A_k(\zeta_1) \partial_j A_i(\zeta_1) (y_i-\zeta'_{1,i}) (y_j-\zeta'_{1,j})(y_k-\zeta'_{1,k}) (y_\ell-\zeta'_{1,\ell})(y_s-\zeta'_{1,s})\F{w'}{|y-\zeta'_1|}\E^{\I\sigma_1+\I \bm A(\zeta_1)\cdot y}  \\
  & - \F{\e^2}{2} \sum_{i,j,k,\ell=1}^N \partial_j A_i(\zeta_1) \partial_{\ell k} A_\ell(\zeta_1) (y_i-\zeta'_{1,i}) (y_j-\zeta'_{1,j}) (y_k-\zeta'_{1,k})  w(|y-\zeta'_1|) \E^{\I\sigma_1+\I \bm A(\zeta_1)\cdot y} \\
  & +\left[ O(\e^3 |y-\zeta'_1|^5 + \e^3 |y-\zeta'_1|^4) + \I O(\e^2 |y-\zeta'_1|^4 )\right] w(|y-\zeta'_1|))\E^{\I\sigma_1+\I \bm A(\zeta_1)\cdot y}.
\end{align*}
Therefore, in $|y-\zeta'_m|\le \frac\delta{\sqrt{\e}}$
\BEN
R(y) \E^{-\I\sigma_m-\I \bm A(\zeta_m)\cdot y}:= R_{m,1}(y) + \I R_{m,2}(y),
\EEN
 where
\begin{align}
& \ R_{m,1}(y) \nonumber\\
=&\ \e^2 \sum_{i,j,k=1}^N \partial_j A_i(\zeta_m) \partial_k A_i(\zeta_m) (y_j-\zeta'_{m,j}) (y_k-\zeta'_{m,k}) w(y-\zeta'_m) \nonumber\\
&\ -\e^2 \sum_{j,k,\ell=1}^N \partial_\ell A_k(\zeta_m) \left[ \partial_j A_k(\zeta_m)+ \partial_k A_j(\zeta_m) \right] (y_j-\zeta'_{m,j})(y_\ell-\zeta'_{m,\ell})w(|y-\zeta'_m|)  \nonumber \\
  &\ - \e^2 \sum_{i,j,k,\ell=1}^N \partial_\ell A_k(\zeta_m)\partial_j A_i(\zeta_m)(y_i-\zeta'_{m,i})(y_j-\zeta'_{m,j})(y_k-\zeta'_{m,k})(y_\ell-\zeta'_{m,\ell})\F{w'(|y-\zeta'_m|)}{|y-\zeta'_m|}  \nonumber\\
  &\ -\F{\e^2}{2}(\nabla_x \cdot \bm A(\zeta_m))\sum_{i,j=1}^N  \partial_j A_i(\zeta_m)(y_i-\zeta'_{m,i})(y_j-\zeta'_{m,j})w(|y-\zeta'_m|) \nonumber\\
  &\ -\F{(p-1)\e^2}{8} \sum_{i,j,k,\ell=1}^N \partial_j A_i(\zeta_m)\partial_\ell A_k(\zeta_m)(y_i-\zeta'_{m,i})(y_j-\zeta'_{m,j}) (y_k-\zeta'_{m,k})(y_\ell-\zeta'_{m,\ell})w^p(|y-\zeta'_m|) \nonumber\\
  &\ +  \e^3  \sum_{i,j,k,\ell=1}^N  \partial_j A_i(\zeta_m) \partial_{k\ell} A_i(\zeta_m)(y_j-\zeta'_{m,j}) (y_k-\zeta'_{m,k})  (y_\ell-\zeta'_{m,\ell})w(|y-\zeta'_m|) \nonumber\\
  &\ - \F{\e^3}{2}  \sum_{j, k,\ell,s=1}^N \partial_{\ell s}A_k(\zeta_m) \left[ \partial_j A_k(\zeta_m) + \partial_k A_i(\zeta_m)\right]  (y_j-\zeta'_{m,j})(y_\ell-\zeta'_{m,\ell})(y_s-\zeta'_{m,s}) w(|y-\zeta'_m|) \nonumber \\
  &\ - \F{\e^3}{2}  \sum_{i,j,k,\ell,s=1}^N  \partial_{\ell s}A_k(\zeta_m) \partial_j A_i(\zeta_m) (y_i-\zeta'_{m,i}) (y_j-\zeta'_{m,j})(y_k-\zeta'_{m,k}) (y_\ell-\zeta'_{m,\ell})(y_s-\zeta'_{m,s})\F{w'(|y-\zeta'_m|)}{|y-\zeta'_m|} \nonumber \\
  &\ - \F{\e^3}{2} \sum_{i,j,k,\ell=1}^N \partial_j A_i(\zeta_m) \partial_{\ell k} A_\ell(\zeta_m) (y_i-\zeta'_{m,i}) (y_j-\zeta'_{m,j}) (y_k-\zeta'_{m,k})  w(|y-\zeta'_m|)  \nonumber \\
  &\ + O\left(\e^4 |y-\zeta'_m|^4  +\e^4 |y-\zeta'_m|^5\right)w(|y-\zeta'_m|) + O (\e^{4} |y-\zeta'_m|^8 w^{p}(|y-\zeta'_m|) )  \label{R1}
\end{align}
and
\[
\begin{split}
  R_{m,2}(y) =&\ \e^2  \sum_{i,j,k=1}^N   \partial_{jk} A_i(\zeta_m)(y_j-\zeta'_{m,j}) (y_k-\zeta'_{m,k})\partial_i w(|y-\zeta'_m|)\\
   & \ + \e^2 \sum_{i,j=1}^N  \partial_{ij}A_i(\zeta_m)(y_j-\zeta'_{m,j})w(|y-\zeta'_m|)\\
  &\ + O\left(\e^3 |y-\zeta'_m|^4  + \e^3 |y-\zeta'_m|^2 \right)w(|y-\zeta'_m|) +O\left( \e^{4} |y-\zeta'_m|^{8}\right)w^p(|y-\zeta'_m|).
\end{split}
\]
Hence we get the  estimate
\[
\sum_{m=1}^K \int_{B(\zeta'_m;\frac\delta{\sqrt{\e}})}|R(y)|^2 dy \le C\e^2.
\]
As for the domain $|y-\zeta'_m|\ge \frac\delta{\sqrt{\e}},  \forall \ m=1,2,\ldots, K$, using the asymptotic behaviour of $w(|y-\zeta'_m|)$, it is easy to see that
\[
\int_{\mathbb{R}^N\setminus \cup_{m=1}^K B(\zeta'_m;\frac\delta{\sqrt{\e}})}|R(y)|^2 dy \le C\E^{-\frac\delta{\sqrt{\e}}}.
\]
The result for $R(y)$ is concluded. \par

As for the estimates of $\partial_{\zeta'_{m,k}} R$ and $\partial_{\sigma_m} R$, one may check it  similarly.
\end{proof}

\section{The linear problem and the nonlinear problem}

This section is devoted to the invertibility of the linear operator $L$ in order to solve problem (\ref{phi}):
\[
L\phi = (\I \nabla + \bm A(\e y))^2 \phi + \phi - (p-1) |\mathcal W|^{p-3}\operatorname{Re}(\overline{\mathcal W}\phi)\mathcal W - |\mathcal W|^{p-1}\phi = - R(y) + N(\phi).
\]
Let $H$ be the Hilbert space as the closure of $C_0^\infty(\mathbb{R}^N, \mathbb{C})$ under the scalar product
\[
( u, v)=\operatorname{Re}\int_{\mathbb{R}^N}\left(\I\nabla u +\bm A(\e y)u\right) \overline{\left(\I\nabla v +\bm A(\e y)v\right)} + u\bar{v}.
\]
The norm deduced by the above scalar product is equivalent to the usual norm of $H^1(\mathbb{R}^N, \mathbb{C})$ due to the boundness of $|\bm A(x)|$, see \cite{CS1}.
In $|y-\zeta'_m|\le \frac\delta{\sqrt{\e}}$, the operator $L$ formally looks like
\BEN
(\I \nabla + \bm A(\zeta_m))^2 \phi + \phi - (p-1)|U_m|^{p-3}\operatorname{Re}(\overline{ U}_m\phi) U_m - | U_m|^{p-1}\phi,
\EEN
which is not invertible. Precisely, the null space of this limit operator is
\BEN
\operatorname{span}_\mathbb{R}\{Z_{m,0}, Z_{m,1},\cdots,Z_{m,N}\}
\EEN
where
\[
Z_{m, 0}=\I w(|y-\zeta_m'|)\E^{\I \sigma_m + \I \bm A(\zeta_m)\cdot y} = \I U_m \qquad \text{and}\qquad Z_{m,i}=\frac{\partial w(|y-\zeta_m'|)}{\partial \zeta'_{m,i}}\E^{\I \sigma_m + \I \bm A(\zeta_m)\cdot y},\quad 1\le i \le N.
\]
The symbol span$_\mathbb{R}$ means  the linear combinations on real numbers, see for instance \cite{CS1,CS2}.
Therefore, we study the following linear problem with  $h \in L^2(\mathbb{R}^N, \mathbb{C})$
\BE\label{linear problem}
\begin{cases}
\D L\phi = h + \sum_{i=0}^N\sum_{m=1}^K c_{m,i} \chi_m Z_{m,i}, \\
\D \operatorname{Re}\int_{\mathbb R^N} \chi_m \overline Z_{m,i} \phi = 0, \quad i=0,1,\cdots,N,\quad m=1, \ldots,  K
\end{cases}
\EE
where $\chi_m(y)=\chi(|y-\zeta'_m|)$ is a smooth  cut-off  function on the large ball $B_R(\zeta'_m)$, satisfying $\chi(s)=1$ for $|s|\le R$ and $\chi(s)=0$ for $|s|\ge R+1$.

%We should point out that the orthogonality to $Z_{m,0}$  is not included in (\ref{linear problem}), since we have fixed the freedom of phase $\E^{\I\sigma}$ to $\E^{\I 0}$ in the original problem (\ref{main problem}).
Next,  we will prove the following  invertibility proposition which is the main result in this section.
\BP\label{t3.1}
The linear problem (\ref{linear problem}) admits a unique solution  $(\phi,c_{m,i})=(\widetilde{T}(h),c_{m,i})$, $i=0,1,\cdots,N,\ m=1, \ldots,  K$  satisfying
\BEN
\|\phi\|_{H^2} = \|\widetilde{T}(h)\| \leq C \|h\|_{L^2}, \qquad |c_{m,i}|\leq C\|h\|_{L^2}.
\EEN
\EP

Before giving the proof, it is necessary to get an apriori estimate.
\BL\label{l1}
If $(\phi,c_{m,i})$ is a solution of  the problem (\ref{linear problem}), then
\BEN
\|\phi\|_{H^2} \leq C \|h\|_{L^2}, \qquad |c_{m,i}|\leq C\|h\|_{L^2} .
\EEN
\EL
%The existence depends on the following two estimates. To the purpose, we study the coming linear problem before (\ref{linear problem})  that
%\BE\label{linear problem1}
%\begin{cases}
%\D L\phi = h + \sum_{i=0}^N\sum_{m=1}^K c_{m,i} \chi_m Z_{m,i}, \\
%\D \operatorname{Re}\int_{\mathbb R^N} \chi_m \overline Z_{m,i} \phi = 0, \quad i=0,\cdots,N,\quad m=1, \ldots,  K.
%\end{cases}
%\EE
%
%
%\BP\label{p2}
%Assume that $(\phi,c_{m,i}), i=0,1,\cdots,N,\quad m=1, \ldots,  K$ satisfies the linear problem (\ref{linear problem1}) for $h \in L^2(\mathbb{R}^N, \mathbb{C})$. Then for small enough $\e$, we have
%\BEN
%\|\phi\|_{H^1} \leq C \|h\|_{L^2}, \qquad |c_{m,i}|\leq C\|h\|_{L^2} + O(\e).
%\EEN
%\EP
\begin{proof}
The proof is very standard and we here prove it briefly for the completion.
First, we test the equation (\ref{linear problem}) by $\overline Z_{\ell, j}$, $1\le \ell \le K$, $0\le j\le N$ and get that
\BE\label{6}
\langle L\phi,  \overline Z_{\ell, j} \rangle= \text{Re}\int_{\mathbb R^N} h \overline Z_{\ell, j} + c_{\ell, j} \int_{\mathbb R^N}  |Z_{\ell,j}|^2+O\left( \E^{-\delta/{\e}} \sum_{i=0}^N\sum_{m=1}^K |c_{m,i}| \right).
\EE
Note that
\[
\begin{split}
\text{Re}\int_{\mathbb R^N} (\I\nabla + A(\e y))\phi \overline{(\I\nabla + A(\e y)) Z_{\ell,j} }
=& \text{Re}\int_{\mathbb R^N}  (\I\nabla +A(\e y))Z_{\ell,j} \overline{(\I\nabla +A(\e y))\phi}\\
=& \text{Re}\int_{\mathbb R^N}  (\I\nabla +A(\zeta_\ell))Z_{\ell,j} \overline{(\I\nabla +A(\zeta_\ell))\phi} + O(\e) \|\phi\|_{H^1},
\end{split}
\]
and
\BEN
\left| \int_{\mathbb R^N} h  Z_{\ell, j} \right| \leq C \|h\|_{L^2}.
\EEN
Thus, it holds from (\ref{6}) and the equation of ${\overline Z_{\ell, j}}$ that
\BE\label{7}
c_{\ell,j} = O(\e \|\phi\|_{H^1} + \|h\|_{L^2}).
\EE

Next we will prove $\|\phi\|_{H^1} \leq C \|h\|_{L^2}$ by contradiction.  Suppose that for some sequence $\{\e_n\}$, there always exist $\phi_n$ and $h_n$ such that
\BEN
\|\phi_n\|_{H^1} = 1 \qquad \text{and} \qquad \|h_n\|_{L^2} =o(1)\ \text{as }\e_n\to 0.
\EEN
Testing  (\ref{linear problem})  against $ \eta_m\varphi (\varphi\in C_c^\infty(\mathbb R^N,\mathbb C))$ where $\eta_m(y)\equiv 1$ in $|y-\zeta_m'|<\F{\delta}{\sqrt{\e}}$ and $\eta_m(y)\equiv 0$ in  $|y-\zeta_m'|>\F{2\delta}{\sqrt{\e}}$, one can obtain that
\begin{eqnarray*}
 && \text{Re}\int_{\mathbb R^N} (\I\nabla + A(\e_n y))\phi_n \overline{(\I\nabla + A(\e_n y)) \eta_m \varphi} + \text{Re}\int_{\mathbb R^N} \phi_n \eta_m  \overline \varphi \\
 &&- (p-1)\text{Re}\int_{\mathbb R^N} |\mathcal W|^{p-3} (\operatorname{Re}(\overline{\mathcal W}\phi_n))\mathcal W \eta_m \overline \varphi - \text{Re}\int_{\mathbb R^N} |\mathcal W|^{p-1} \phi_n \eta_m \overline \varphi  \\
  = &&\text{Re}\int_{\mathbb R^N} h_n \eta_m \overline \varphi +  \sum_{i=0}^N c_{m,i} \text{Re}\int_{\mathbb R^N} \chi_m Z_{m, i} \overline \varphi.
\end{eqnarray*}
Note that $\phi_n \rightharpoonup \phi$ in $H^1_{\text{loc}}(\mathbb{R}^N, \mathbb{C})$ up to a subsequence. Thus dominated convergence theorem tells us that
\begin{eqnarray*}
&& \text{Re}\int_{\mathbb R^N} (\I\nabla + A(\zeta_m))\phi \overline{(\I\nabla + A(\zeta_m))  \varphi} + \text{Re}\int_{\mathbb R^N} \phi \bar \varphi  \\&& - (p-1)\text{Re}\int_{\mathbb R^N} |U_m|^{p-3}(\operatorname{Re}(\overline{ U}_m\phi)) U_m  \bar \varphi - \text{Re}\int_{\mathbb R^N} |U_m|^{p-1} \phi  \bar \varphi = 0.
\end{eqnarray*}
This means that  $\phi$ is a  solution of
\BEN
(\I \nabla + \bm A(\zeta_m))^2 \phi +  \phi - (p-1)|U_m|^{p-3}\text{Re}(\overline{U}_m\phi) U_m - |U_m|^{p-1} \phi = 0 \qquad \text{in }\mathbb R^N.
\EEN
Then one gets that $\phi =0$ from the orthogonal conditions, which further implies that
\BE\label{8}
\phi_n \to 0 \quad \text{a.e. in }B_R(\zeta'_m), \qquad \forall\ R>0,\ m=1,2,\ldots,K.
\EE
On the other hand, note  that
\begin{eqnarray}\label{9}
&&\int_{\mathbb R^N} |(\I \nabla + \bm A(\e y))\phi_n|^2 + \int_{\mathbb R^N}|\phi_n|^2 - (p-1)\int_{\mathbb R^N} |\mathcal W|^{p-3}(\operatorname{Re}(\overline{\mathcal W}\phi_n))^2 -\int_{\mathbb R^N} |\mathcal W|^{p-1} |\phi_n|^2\nonumber\\
 =&& \operatorname{Re} \int_{\mathbb R^N} \bar h_n \phi_n = o(1).
\end{eqnarray}
From (\ref{8}) and the exponential decay of $|U_m|$, we obviously have
\BEN
\int_{\mathbb R^N} |\mathcal W|^{p-1} |\phi_n|^2 = \sum_{m=1}^N\int_{B_R(\zeta'_m)} |U_m|^{p-1} |\phi_n|^2 + \int_{\mathbb R^N\setminus \cup_{m=1}^K B_R(\zeta'_m)} |\mathcal W|^{p-1} |\phi_n|^2+O(\e)=O(\E^{-R})+o(1).
\EEN
So is $\int_{\mathbb R^N} |\mathcal W|^{p-3}(\operatorname{Re}(\overline{\mathcal W}\phi_n))^2$, which together with (\ref{9}) shows that
\BEN
\int_{\mathbb R^N}|\phi_n|^2 = O(e^{-R})+o(1) \quad\text{and}\quad \int_{\mathbb R^N} |(\I \nabla + \bm A(\e y))\phi_n|^2 = O(\E^{-R})+o(1).
\EEN
Finally, it is derived that
\begin{eqnarray*}
O(\E^{-R})+o(1) &=& \int_{\mathbb R^N} |(\I \nabla + \bm A(\e y))\phi_n|^2 \\
&=&  \int_{\mathbb R^N} |\nabla \phi_n|^2 + \int_{\mathbb R^N} |\bm A(\e_n y)|^2|\phi_n|^2 + 2\text{Re}\int_{\mathbb R^2} \I \bm A(\e_n y) \cdot \nabla\phi_n \phi_n \\
&\geq& \F{1}{2}  \int_{\mathbb R^N} |\nabla \phi_n|^2 -  \int_{\mathbb R^N} |\bm A(\e_n y)|^2|\phi_n|^2 = \F{1}{2} \|\phi_n\|_{H^1} + O(\E^{-R})+o(1),
\end{eqnarray*}
since  $\bm A(x)$ is bounded. This leads to a  contradiction to $\|\phi_n\|_{H^1} =1$. Hence
$$\|\phi\|_{H^1} \leq C\|h\|_{L^2}.$$

At last, the regularity theory gives  $\|\phi\|_{H^2} \leq C\|h\|_{L^2}$.
\end{proof}

\begin{proof}[Proof of Proposition \ref{t3.1}]
Denote the Hilbert space
\BEN
\mathcal H = \left\{ \phi \in H(\mathbb R^N) ~\left|~ \text{Re}\int_{\mathbb R^N} \chi_m \phi \overline Z_{m,i} = 0 \right. \right\}.
\EEN
Then, from  Riesz representation theorem the equation (\ref{linear problem}) is equivalent to
\BEN
\phi - T(\phi) = \tilde h, \qquad \text{in }\mathcal H,
\EEN
where $T$ is a compact operator on $\mathcal H$.
Based on Proposition \ref{t3.1}, Fredholm alternative tells us the unique existence of $\phi$. And $c_{m,i}$ can be given by  $\phi$ using integration. Their estimates were given in the above proposition.
\end{proof}

Also for $\phi=\widetilde{T}(h)$, it is important for later purposes to understand the differentiability of the operator  $\widetilde{T}$ with respect to  $\zeta_j'$ and $\sigma_j, j=1, \ldots, K$. Recall that $\phi$ satisfies the equation
\BEN
L\phi= (\I\nabla + A(\e y))^2\phi + \phi -(p-1)|\mathcal W|^{p-3}\text{Re}(\overline{\mathcal W} \phi) \mathcal W- |\mathcal W|^{p-1}\phi = h + \sum_{m,i} c_{m,i} \chi_m Z_{m,i}.
\EEN
Thus
\begin{eqnarray*}
L(\partial_{\zeta'_{j,k}}\phi) &=& (\I\nabla + A(\e y))^2(\partial_{\zeta'_{j,k}}\phi) + \partial_{\zeta'_{j,k}}\phi -(p-1)|\mathcal W|^{p-3}\text{Re}(\overline{\mathcal W} \partial_{\zeta'_{j,k}}\phi) \mathcal W - |\mathcal W|^{p-1} \partial_{\zeta'_{j,k}}\phi \\
&=& O(|\mathcal W|^{p-2}  |\phi| |\partial_{\zeta'_{j,k}} \mathcal W| ) + \sum_{i} c_{j,i} \partial_{\zeta'_{j,k}} (\chi_j Z_{j,i}) + \sum_{m,i} (\partial_{\zeta'_{j,k}} c_{m,i})\chi_m Z_{m,i}.
\end{eqnarray*}
Moreover, the derivative of the orthogonal condition is
\begin{gather*}
\text{Re}\int_{\mathbb R^N} \chi_m \overline Z_{m,i} (\partial_{\zeta'_{j,k}}\phi) = 0, \qquad \text{for }m\neq j,  \\
\text{Re}\int_{\mathbb R^N} \chi_j \overline Z_{j,i} (\partial_{\zeta'_{j,k}}\phi) =  - \text{Re}\int_{\mathbb R^N} \partial_{\zeta'_{j,k}} (\chi_j \overline Z_{j,i}) \phi.
\end{gather*}
Set $\varphi= \partial_{\zeta'_{j,k}} \phi +\sum_{i=1}^N b_{jk, i}  \chi_j Z_{j,i}$ and
\BEN
b_{jk,i} = \text{Re}\int_{\mathbb R^N} \partial_{\zeta'_{j,k}} (\chi_j \overline Z_{j,i}) \phi \Big/ \int_{\mathbb R^N} \chi_j |Z_{j,i}|^2.
\EEN
Note that $|b_{jk, i}| \leq C \|\phi\|_{L^2}\leq C\|h\|_{L^2}$.
 Then $\varphi$ satisfies all the orthogonal conditions and
\BEN
L\varphi = b_{ji} L (\chi_j Z_{j,i}) + O(|\mathcal W|^{p-2}  |\phi| |\partial_{\zeta'_{j,k}} \mathcal W| ) + \sum_{i} c_{j,i} \partial_{\zeta'_{j,k}} (\chi_j Z_{j,i}) + \sum_{m,i} (\partial_{\zeta'_{j,k}} c_{m,i})\chi_m Z_{m,i}.
\EEN
With Lemma \ref{l1} in hand,
\BEN
\|\varphi\|_{H^2} \leq C \|b_{ji} L (\chi_j Z_{j,i})\|_{L^2} + C\| |\mathcal W|^{p-2}  \phi \partial_{\zeta'_{j,k}} \mathcal W \|_{L^2} + C \sum_i \|c_{j,i} \partial_{\zeta'_{j,k}} (\chi_j Z_{j,i})\|_{L^2} \leq C \|h\|_{L^2}.
\EEN
Therefore, we conclude that
\BEN
\|\partial_{\zeta'_{j,k}} \widetilde{T} (h)\|_{H^2} = \|\partial_{\zeta'_{j,k}} \phi\|_{H^2} \leq \|\varphi\|_{H^2} + \|b_{ji}  \chi_j Z_{j,i}\|_{H^2} \leq C\|h\|_{L^2}.
\EEN
The same process may be carried out for $\partial_{\sigma_j}\phi$. Based on the above discussion,   the following proposition holds obviously.

\BP\label{p3.3}
For the unique solution $\phi=\widetilde{T}(h)$ in Proposition \ref{t3.1}, it holds that
\BEN
\|\partial_{\zeta'_{j,k}} \widetilde{T}(h)\|_{H^2} \leq C\|h\|_{L^2}, \qquad  \|\partial_{\sigma_j} \widetilde{T} (h)\|_{H^2} \leq C\|h\|_{L^2}, \quad \forall\ j=1,\ldots, K, k=1,\ldots, N.
\EEN
\EP

Now we  can deal with the following nonlinear problem
\BE\label{nonlinear problem}
\begin{cases}
\D L\phi = -R(y) + N(\phi) + \sum_{i=0}^N \sum_{m=1}^K c_{m,i} \chi_m Z_{m, i}, \\
\D \operatorname{Re}\int_{\mathbb R^N} \chi_m \overline Z_{m,i} \phi = 0, \quad i=0,\cdots, N, m=1, \ldots, K.
\end{cases}
\EE

\BP\label{p3.2}
The nonlinear problem (\ref{nonlinear problem}) admits a unique solution $\phi$ satisfying
\BEN
\|\phi\|_{H^2}  = O(\e^{2}).
\EEN
Moreover, $( \bm \sigma, \bm\zeta') \to\phi$ is of class $C^1$ for $ \bm\sigma=(\sigma_1, \ldots, \sigma_K), \bm \zeta'=(\zeta_1', \ldots, \zeta_K')$, and
\BEN
\|\partial_{\zeta'_{j,k}}\phi\|_{H^2}  = O(\e^{2\beta}) \quad \text{and}\quad \|\partial_{\sigma_j}\phi\|_{H^2}  = O(\e^{2\beta}), \qquad \forall\ j,k.
\EEN
\EP

\begin{proof}

Recall that $\beta=\min\{p-1, 1\}$.  The proof is based on the contraction mapping theorem. First, for a large enough number $\gamma_0>0$,  we set
\BEN
\mathcal S = \left\{ \phi \in \mathcal H ~\left|~  \|\phi\|_{H^2} \leq \gamma_0 \e^2\right.\right\}.
\EEN
The  nonlinear problem (\ref{nonlinear problem}) is transferred  to solving
\BEN
\phi = \widetilde{T}(-R(y) + N(\phi)) := \mathcal A(\phi),
\EEN
which means to find a fixed point of the operator $\mathcal A$.\par

First, the operator $\mathcal A$ is from $\mathcal S$ to itself. In fact,
\BEN
\|\mathcal A(\phi)\|_{H^2} = \|\widetilde{T}(-R(y) + N(\phi))\|_{H^2} \leq C\|R(y)\|_{L^2} + C\|N(\phi)\|_{L^2}\leq C\e^2 + C\|\phi\|_{H^2}^{1+\beta} \leq \gamma_0 \e^2.
\EEN
Next the operator $\mathcal A$ is a contraction mapping, since
\begin{eqnarray*}
\|\mathcal A(\phi_1) - \mathcal A(\phi_2)\|_{H^2} &=& \|\widetilde{T} (N(\phi_1)-N(\phi_2))\|_{H^2}\leq C\|N(\phi_1)-N(\phi_2)\|_{L^2}\\
 &\leq& C\left(\|\phi_1\|_{H^2}^\beta + \|\phi_2\|_{H^2}^\beta\right) \|\phi_1-\phi_2\|_{H^2}.
\end{eqnarray*}
Thus $\mathcal A$ has a unique fixed point in $\mathcal S$, which is the unique solution of problem (\ref{nonlinear problem}).\par

Next we come to $\partial_{\zeta'_{j,k}} \phi$ and $\partial_{\sigma_j} \phi$. The $C^1$-regularity in $\zeta'_j$ and $\sigma_j$ is guaranteed by the implicit function theorem. One may refer to the proof of Lemma 4.1 in \cite{CS1}. It is easy to see
\BE\label{12}
\partial_{\zeta'_{j,k}} \phi = \partial_{\zeta'_{j,k}} \widetilde{T} (-R(y) + N(\phi)) + \widetilde{T}(- \partial_{\zeta'_{j,k}} R(y) +\partial_{\zeta'_{j,k}} N(\phi) ).
\EE
Notice  that $\left|\partial_{\zeta'_{j,k}} N(\phi)\right| = O\left( |\phi|^\beta |\partial'_{jk} w|+|\phi|^\beta |\partial_{\zeta'_{j,k}} \phi|\right)$.  So we get
\BEN
\|\partial_{\zeta'_{j,k}} N(\phi)\|_{L^2} \leq C\|\phi\|_{H^2}^{\beta} + C\|\phi\|_{H^2}^{\beta} \|\partial_{\zeta'_{j,k}} \phi\|_{H^2}.
\EEN
Thus (\ref{12}) and Proposition \ref{p3.3} lead to
\BEN
\|\partial_{\zeta'_{j,k}} \phi\|_{H^2} \leq C\|\partial_{\zeta'_{j,k}} R(y)\|_{L^2} + C\|R(y)\|_{L^2} + C\|N(\phi)\|_{L^2} + C\|\partial_{\zeta'_{j,k}} N(\phi)\|_{L^2} \leq C\e^{2\beta}.
\EEN
The  estimate  for $\|\partial_{\sigma_j} \phi\|_{H^2}$ may be gotten by the same process.
\end{proof}

\section{Variational Reduction}
 According to the above discussion, the remaining thing is to let $c_{m,i}=0$ in the nonlinear problem (\ref{nonlinear problem}) in order to make $\mathcal W +\phi$ be a solution of the original problem. It can be done by the variational reduction process.\par

Note the energy functional of problem (\ref{appro}) is
\[
E(u) = \F{1}{2}\int_{\mathbb R^N} |\I \nabla u+ \bm A(\e y) u|^2 \mathrm d y + \F{1}{2}\int_{\mathbb R^N} |u|^2 - \F{1}{p+1}\int_{\mathbb R^N} |u|^{p+1}.
\]
Define
\[
F(\bm \sigma, \bm \zeta')= E(\mathcal W + \phi)(\bm\sigma,\bm\zeta'),
\]
then the existence of critical points to $E(u)$ may be reduced to find critical points of the finite dimensional function $F(\bm \sigma, \bm \zeta')$.

\BP\label{p1}
If $(\bm \sigma,\bm\zeta')$ is a critical point of $F(\bm \sigma, \bm \zeta')$, then
$c_{m,i}=0$ for all $m,i$.
\EP

\begin{proof}
It is easy to see that
\begin{align*}
  \partial_{\zeta_{m,j}'} F(\bm \sigma,\bm \zeta') &=   \partial_{\zeta'_{m,j}} E(\mathcal W + \phi)=  E'(\mathcal W + \phi)\left[\F{\partial \mathcal W}{\partial{\zeta'_{m,j}}}+\F{\partial \phi}{\partial{\zeta'_{m,j}}}\right]\\
  &= \sum_{i,\ell} \text{Re}\int_{\mathbb R^N}c_{\ell,i} \chi_\ell \overline Z_{\ell,i} \left[\F{\partial U_m}{\partial{\zeta'_{m,j}}} + \e \F{\partial \Psi_m}{\partial{\zeta'_{m,j}}} + \F{\partial \phi}{\partial{\zeta'_{m,j}}} \right] =  c_{m,j} \int_{\mathbb R^N}\chi_m |Z_{m,j}|^2 +o(1).
\end{align*}
Similarly, it is also true that
\BEN
\partial_{\sigma_m} F(\bm \sigma,\bm \zeta') = c_{m,0} \int_{\mathbb R^N}\chi_m |Z_{m,0}|^2 +o(1).
\EEN
Thus $c_{m,i} = 0$ if $(\bm\sigma, \bm \zeta')$ is a critical point of $F$ since the coefficient matrix of $c_{m,i}$ is diagonal dominant.
\end{proof}

Next we should calculate $F(\bm\sigma,\bm\zeta')$ in view of Proposition \ref{p1}.
\BP\label{p4.2}
It holds that
\begin{gather*}
F(\bm \sigma, \bm \zeta') = E(\mathcal W) +O(\e^4),
\intertext{and for any $j,k$,}
\partial_{\zeta'_{j,k}} F(\bm \sigma, \bm \zeta') = \partial_{\zeta'_{j,k}} E(\mathcal W) + O(\e^{2+2\beta}), \qquad \partial_{\sigma_j} F(\bm \sigma, \bm \zeta') = \partial_{\sigma_j} E(\mathcal W) + O(\e^{2+2\beta}).
\end{gather*}
\EP

\begin{proof}
Direct computation leads to
 \BEN
\begin{split}
F(\bm \sigma, \bm \zeta') =& \F{1}{2}\int_{\mathbb R^N} \left|\I \nabla (\mathcal W+\phi) + \bm A(\e y) (\mathcal W+ \phi)\right|^2 + \F{1}{2}\int_{\mathbb R^N} |\mathcal W + \phi|^2 - \F{1}{p+1}\int_{\mathbb R^N} |\mathcal W + \phi|^{p+1}  \\
=&\ E(\mathcal W) + \F{1}{2}\int_{\mathbb R^N}\text{Re} \left((R(y)-N(\phi)) \bar\phi\right)  -\int_{\mathbb R^N} \Big[\F{1}{p+1}|\mathcal W+\phi|^{p+1} - \F{1}{p+1}|\mathcal W|^{p+1} \\
& - |\mathcal W|^{p-1} \text{Re}(\mathcal W \bar \phi) - \F{p-1}{2} |\mathcal W|^{p-3}(\text{Re}(\mathcal W \bar \phi))^2 - \F{1}{2} |\mathcal W|^{p-1}|\phi|^2\Big].
\end{split}
\EEN
Then the proposition follows from Proposition \ref{p2.1} and Proposition \ref{p3.2} easily. By the computation in Proposition \ref{p3.2}, it is easy to check that
\BEN
\partial_{\zeta'_{j,k}} F(\bm \sigma, \bm \zeta') = \partial_{\zeta'_{j,k}} E(\mathcal W) + O(\|\phi\|^{1+\beta}).
\EEN
So is $\partial_{\sigma_j} F(\bm \sigma, \bm \zeta')$.
\end{proof}

Since $E(\mathcal W)$ is the main part of $F(\bm\sigma, \bm\zeta')$, it is important to get the expression of $E(\mathcal W)$. Elegant computation shows the following proposition.

\BP\label{p4.3}
It holds that for $\e $ small enough,
\BEN
E(\mathcal W)=A_0K + B_0 \e^2 \sum_{m=1}^K \sum_{i,j=1}^N (\partial_i A_j(\zeta_m) - \partial_j A_i(\zeta_m))^2+ O(\e^{4}).
\EEN
Furthermore, the remainder term $O(\e^{4})$ also holds for the derivatives in $\bm\zeta', \bm\sigma$.
Here $A_0 = \F{p-1}{2(p+1)} \int_{\mathbb R^N} w^{p+1}(|y|) \mathrm dy$ and $B_0 = \F{1}{4} \int_{\mathbb R^N} y_1^2 w^2(|y|) \mathrm dy$ are both universal positive constants.
\EP

\begin{proof}
It is easy to see that
\begin{align}
E(\mathcal W)  &= \F{1}{2}\int_{\mathbb R^N} |\I \nabla \mathcal W+ \bm A(\e y)\mathcal W|^2 + \F{1}{2}\int_{\mathbb R^N} |\mathcal W|^2 - \F{1}{p+1}\int_{\mathbb R^N} |\mathcal W|^{p+1}\nonumber\\
&= \F{1}{2} \operatorname{Re} \int_{\mathbb R^N} R(y) \overline{\mathcal W} + \F{p-1}{2(p+1)} \int_{\mathbb R^N} |\mathcal W|^{p+1} \nonumber\\
&=\sum_{m=1}^K \text{Re}\int_{B(\zeta'_m;\frac\delta{\sqrt{\e}})}\left[\frac12 R(y) \overline{\mathcal W} + \frac{p-1}{2(p+1)} |\mathcal W|^{p+1}\right]+ \text{Re}\int_{\mathbb R^N\setminus \cup_{m=1}^K B(\zeta'_m;\frac\delta{\sqrt{\e}}) }\left[\frac12 R(y) \overline{\mathcal W} + \frac{p-1}{2(p+1)} |\mathcal W|^{p+1}\right] \nonumber\\
&= \sum_{m=1}^K  \int _{B(\zeta'_m;\frac\delta{\sqrt{\e}})} \Bigg[\frac{ R_{m,1}}2 w(|y-\zeta'_m|) + \e \frac{R_{m,2}}2 \psi_m(y) +\frac{p-1}{2(p+1)} |U_m|^{p+1}\nonumber\\
&\qquad + \F{(p-1)\e^2}{4}  |U_m|^{p-1} |\Psi_m|^2\Bigg]  + O(\e^{4} ),\label{10}
\end{align}
where $\psi_m, \Psi_m$ are given in (\ref{psi}), and $R_{m,1}$, $R_{m,2}$ are defined in Proposition \ref{p2.1}.
First, by the oddness of the terms in order $\e^3$ in (\ref{R1}), it can be  obtained  that
\begin{align*}
  &\ \int _{B(\zeta'_m; \frac\delta{\sqrt{\e}})} R_{m,1} w(|y-\zeta'_m|) \nonumber \\
  =&\ \e^2 \sum_{i,j,k=1}^N \partial_j A_i(\zeta_m) \partial_k A_i(\zeta_m) \int_{B(\zeta'_m;\frac\delta{\sqrt{\e}})} (y_j-\zeta'_{m,j}) (y_k-\zeta'_{m,k}) w^2(|y-\zeta'_m|) \nonumber \\
  &\ -\e^2\sum_{j,k,\ell=1}^N \partial_\ell A_k(\zeta_m)\left [\partial_j A_k(\zeta_m)+ \partial_k A_j(\zeta_m)\right] \int_{B(\zeta'_m;\frac\delta{\sqrt{\e}})} (y_j-\zeta'_{m,j})(y_\ell-\zeta'_{m,\ell})w^2(|y-\zeta'_m|) \nonumber \\
  &\ -\e^2 \sum_{i,j,k,\ell}^N\partial_\ell A_k(\zeta_m)\partial_j A_i(\zeta_m)\int_{B(\zeta'_m;\frac\delta{\sqrt{\e}})}(y_i-\zeta'_{m,i})(y_j-\zeta'_{m,j})(y_k-\zeta'_{m,k})(y_\ell-\zeta'_{m,\ell})\nonumber \\
   &\ \qquad \qquad \qquad \qquad \qquad \qquad \qquad \qquad \qquad \cdot\F{w'(|y-\zeta'_m|)}{|y-\zeta'_m|}w(|y-\zeta'_m|)  \nonumber \\
  &\ -\F{\e^2}{2}(\nabla_x \cdot \bm A(\zeta_m))\sum_{i,j=1}^N  \partial_j A_i(\zeta_m)\int_{B(\zeta'_m;\frac\delta{\sqrt{\e}})}(y_i-\zeta'_{m,i})(y_j-\zeta'_{m,j})w^2(|y-\zeta'_m|)  \nonumber \\
  &\ -\F{(p-1)\e^2}{8} \sum_{i,j,k,\ell}^N  \partial_j A_i(\zeta_m)\partial_\ell A_k(\zeta_m)\int_{B(\zeta'_m;\frac\delta{\sqrt{\e}})} (y_i-\zeta'_{m,i})(y_j-\zeta'_{m,j}) (y_k-\zeta'_{m,k})(y_\ell-\zeta'_{m,\ell})w^{p+1}  \nonumber \\
  &\  + O(\e^{4}).
%  =&\ \e^2 \sum_{i,j=1}^N (\partial_j A_i(\zeta_m))^2 \int_{\mathbb R^N} (y_j-\zeta'_{m,j})^2 w^2(|y-\zeta'_m|)  \nonumber \\
%  &\ - \e^2 \sum_{i,j=1}^N \partial_j A_i (\zeta_m) (\partial_j A_i (\zeta_m) + \partial_i A_j(\zeta_m)) \int_{\mathbb R^N} (y_j-\zeta'_{m,j})^2 w^2(|y-\zeta'_m|)  \nonumber \\
%  &\ -\e^2 \sum_{i,j=1}^N \partial_j A_i(\zeta_m)(\partial_j A_i(\zeta_m)+\partial_i A_j(\zeta_m))\int_{\mathbb R^N}(y_i-\zeta'_{m,i})^2(y_j-\zeta'_{m,j})^2\F{w'(|y-\zeta'_m|)}{|y-\zeta'_m|}w(|y-\zeta'_m|)  \nonumber \\
%  &\quad -\F{\e^2}{2}(\nabla_x \cdot \bm A(\zeta_m))\sum_{i=1}^N  \partial_i A_i(\zeta_m)\int_{\mathbb R^N}(y_i-\zeta'_{m,i})^2 w^2(|y-\zeta'_m|) \nonumber \\
%  &\ -\F{(p-1)\e^2}{8}\sum_{i,j=1}^N \partial_j A_i(\zeta_m) (\partial_j A_i(\zeta_m)+ \partial_i A_j(\zeta_m) ) \int_{\mathbb R^N} (y_i-\zeta'_{m,i})^2 (y_j-\zeta'_{m,j})^2 w^{p+1}(|y-\zeta'|) \nonumber \\
%  &\ +O(\e^{4}).
\end{align*}

Since integration by parts gives
\begin{eqnarray*}
\int_{\mathbb R^N} y_i y_j y_k y_\ell w(|y|) \F{w'(|y|)}{|y|} \mathrm dy &=& \F{1}{2} \int_{\mathbb R^N} y_i y_j y_k  (\partial_{\ell} w^2) \\
&=& - \F{1}{2} \int_{\mathbb R^N} \delta_{i\ell} \delta_{jk} y_j y_k w^2 - \F{1}{2} \int_{\mathbb R^N} \delta_{j\ell} \delta_{ik} y_i y_k w^2 - \F{1}{2} \int_{\mathbb R^N} \delta_{k\ell} \delta_{ij} y_i y_j w^2 ,
\end{eqnarray*}
%\BEN
%\int_{\mathbb R^N} y_i^2 y_j^2 w(|y|) \F{w'(|y|)}{|y|} \mathrm dy= \int_{\mathbb R^N} y_i^2 y_j w \partial_j w \mathrm dy= \F{1}{2}\int_{\mathbb R^N} y_i^2 y_j \partial_j w^2 \mathrm dy= - \F{1}{2}\int_{\mathbb R^N} \left[y_i^2+2\delta_{ij}y_iy_j\right]  w^2\mathrm dy,
%\EEN
one gets
\begin{eqnarray*}
  &&\sum_{i,j,k,\ell=1}^N\partial_\ell A_k(\zeta_m)\partial_j A_i(\zeta_m) \int_{\mathbb R^N}(y_i-\zeta'_{m,i})(y_j-\zeta'_{m,j})(y_k-\zeta'_{m,k})(y_\ell-\zeta'_{m,\ell})
  \F{w'(|y-\zeta'_m|)}{|y-\zeta'_m|}w(|y-\zeta'_m|) \\
  =&& - 2 B_0\sum_{i,j} \left[\partial_i A_i(\zeta_m) \partial_j A_j(\zeta_m) + (\partial_j A_i(\zeta_m))^2 +  \partial_i A_j(\zeta_m) \partial_j A_i(\zeta_m) \right].
\end{eqnarray*}

It may be checked that
\begin{align}
 &\ \int _{B(\zeta'_m;\frac\delta{\sqrt{\e}})} \F{R_{m,1}}{2} w(|y-\zeta'_m|)\nonumber \\
 =&\ B_0 \e^2 \Bigg \{ 2\sum_{i,j}(\partial_j A_i(\zeta_m))^2- 2\sum_{i,j} \partial_j A_i(\zeta_m) (\partial_j A_i(\zeta_m) + \partial_i A_j(\zeta_m)) \nonumber \\
  &\qquad\quad   + \sum_{i,j} \left[\partial_i A_i(\zeta_m) \partial_j A_j(\zeta_m) + (\partial_j A_i(\zeta_m))^2 +  \partial_i A_j(\zeta_m) \partial_j A_i(\zeta_m) \right] -  \left( \sum_i \partial_i A_i (\zeta_m)\right)^2 \Bigg\}\nonumber \\
&\ -\F{(p-1)\e^2}{16} \sum_{i,j,k,\ell}^N  \partial_j A_i(\zeta_m)\partial_\ell A_k(\zeta_m)\int_{B(\zeta'_m;\frac\delta{\sqrt{\e}})} (y_i-\zeta'_{m,i})(y_j-\zeta'_{m,j}) (y_k-\zeta'_{m,k})(y_\ell-\zeta'_{m,\ell})w^{p+1} \nonumber\\
  &\ +O(\e^{4}). \label{4}
\end{align}
Also
\begin{multline*}
   \int _{B(\zeta'_m;\frac\delta{\sqrt{\e}})} |U_m|^{p-1}|\Psi_m|^2 \\
    = \F{1}{4} \sum_{i,j,k,\ell=1}^N  \partial_j A_i(\zeta_m)\partial_\ell A_k(\zeta_m) \int _{B(\zeta'_m;\delta)} (y_i-\zeta'_{m,i})(y_j-\zeta'_{m,j})(y_k-\zeta'_{m,k})(y_\ell-\zeta'_{m,\ell})w^{p+1}(|y-\zeta'_m|)  + O(\E^{-\frac\delta\e}),
\end{multline*}
and
\BEN
\int_{B(\zeta'_m;\frac\delta{\sqrt{\e}})} R_{m,2} \psi_m  =   O(\e^{3})
\EEN
by the oddness.
Obviously one  has
 $$\int_{B(\zeta'_m;\frac\delta{\sqrt{\e}}) } |U_m|^{p+1} = \int_{\mathbb R^N} w^{p+1}(y) \mathrm dy + O(\E^{-\frac\delta\e}).$$  Note that in (\ref{4})  the term containing $w^{p+1}$ is canceled with $\F{p-1}{4}\e^2\int _{\mathbb R^N} |U_m|^{p-1}|\Psi_m|^2$.
So we conclude, from (\ref{10}), that
\begin{align*}
E(\mathcal W) &= A_0K + \e^2 B_0 \sum_{m=1}^K \sum_{i,j=1}^N \left[ (\partial_j A_i(\zeta_m))^2 - \partial_j A_i(\zeta_m) \partial_i A_j (\zeta_m) \right]+ O(\e^{4}) \\
&= A_0K + \e^2 \F{B_0}{2}  \sum_{m=1}^K \sum_{i, j=1 }^N \left ((\partial_j A_i(\zeta_m)) - \partial_i A_j(\zeta_m)\right )^2 + O(\e^{4}).
\end{align*}
The last equality is due to the symmetry of the indexes $i$ and $j$.
The remainder $O(\e^4)$ also holds for the derivatives of $E(\mathcal W)$ in $(\bm \sigma, \bm\zeta')$ from directly checking the expressions of $R_{m,1}$ and $R_{m,2}$ in the proof of Proposition \ref{p2.1}.

\end{proof}

\section{Proof of the main theorems}
This section devotes to the proof of main theorems.\par

\begin{proof}[Proof of Theorem \ref{t1}]
  Proposition \ref{p4.2} and Proposition \ref{p4.3} mean
\BEN
F(\bm \sigma, \bm\zeta') = A_0K + B_0 \e^2\sum_{m=1}^K \|\bm B(\zeta_m)\|_F^2 + O(\e^{4}).
\EEN
We shall show that $F$ has a critical point under the assumption. Note that for any fixed $\bm\zeta'$, $F(\bm\sigma,\bm\zeta')$ is periodic in $\bm\sigma\in([0,2\pi])^K$. So there always exists a $\bm\sigma(\bm\zeta')$ such that $\partial_{\bm\sigma} F(\bm\sigma(\bm\zeta'),\bm\zeta')=0$.
Next consider the configuration set $\bm\Omega' =\Omega'_1\times\Omega'_2\times\cdots\times\Omega'_m$ of $\bm\zeta'=(\zeta'_1,\zeta'_2,\cdots,\zeta'_m)$, where $\Omega'_m = \e^{-1}\Omega_m$. Obviously
\BEN
\max_{\bm \Omega'} F(\bm\sigma(\bm\zeta'),\bm\zeta') \geq A_0 K + B_0 \e^2 \sum_{m=1}^K \|\bm B(P_m)\|_F^2 + O(\e^4).
\EEN
On the other hand, for any $\bm\zeta'$ on the boundary $\partial\bm\Omega'$, i.e. at least $\zeta_{1}'\in\partial\Omega_{1}'$ without loss of generality, then $\|B(\zeta_1)\|^2_F \leq \|B(P_1)\|^2_F -\delta_0$ for some fixed small $\delta_0>0$. Thus one finds that
\BEN
F(\bm\sigma(\bm\zeta'),\bm\zeta')\Big|_{\bm\zeta'\in\partial\bm\Omega'} \leq  A_0 K + B_0 \e^2 \sum_{m=2}^K \|\bm B(P_m)\|_F^2 + B_0 \e^2 (\|\bm B(P_1)\|_F^2-\delta_0) +O(\e^4).
\EEN
Therefore  $\max_{\bm \Omega'} F(\bm\sigma(\bm\zeta'),\bm\zeta') > \max_{\partial\bm \Omega'} F(\bm\sigma(\bm\zeta'),\bm\zeta')  $.
It implies that $F(\bm\sigma,\bm\zeta')$ admits a critical point.\par

The same procedure can   be carried out  for the case of $K$ local minimum points. Theorem \ref{t1} concludes from  Proposition \ref{p1}.
\end{proof}

\begin{proof}[Proof of Theorem \ref{mt}]
From Proposition \ref{p4.2} and Proposition \ref{p4.3}, we see that
\BEN
\nabla_{\zeta_m'} F(\bm\sigma,\bm\zeta') = B_0 \e^2 \nabla_{\zeta_m'} (\|\bm B(\zeta_m)\|_F^2) + O(\e^{2+2\beta})= B_0 \e^3 \nabla_{\zeta_m} (\|\bm B(\zeta_m)\|_F^2) + O(\e^{2+2\beta}) .
\EEN
Assume $m=1$ for simplicity.  Choose $\zeta_1 = P_1 + \e^\alpha \xi_1$  where $0<\alpha<2\beta-1$. Here the assumption $p>\frac32$  is used to let $\beta=\min\{p-1,  1\}>1/2$. Then it is equivalent to find a  $|\xi_1|\leq 1$ such that
\begin{eqnarray*}
0&=&\nabla_{\zeta_1} (\|\bm B(\zeta_1)\|_F^2)+ O(\e^{2\beta-1}) \\
&=& \nabla_{\zeta_1} (\|\bm B(P_1)\|_F^2) + \e^\alpha \nabla^2_{\zeta_1\zeta_1} (\|\bm B(P_1)\|_F^2) \cdot \xi_1 + O(\e^{2\alpha}|\xi_1|^2) +O(\e^{2\beta-1}).
\end{eqnarray*}
Thus, the nondegeneracy of the critical point $P_1$ leads to the existence of $\xi_1, |\xi_1|=o(1)$ from the Brouwer fixed point theorem.
 Lastly, the existence of critical  $\bm\sigma$ is guaranteed by the periodicity just like in the proof of Theorem \ref{t1}. The proof is complete.
 %Since $F(\bm\sigma, \bm\zeta)$ is periodic  in $[0, 2\pi]^K$ and $C^1$, we just focus on $\nabla_{\bm\zeta} F(\bm\sigma, \bm\zeta)$. Then the following can be obtained\underline{}
%\BEN
%\nabla_{\bm\zeta} F(\bm\sigma, \bm\zeta) = \e^2 \nabla_{\bm\zeta}  \Gamma(\bm\zeta) + O(\e^4),
%\EEN
%where $\Gamma(\bm \zeta)$ is given by (\ref{gamma}). Using the assumption of $\bm\xi$ with the notation $\bm\zeta=\bm\xi+\e\tilde{\bm\zeta}, \ \bm\xi=(\xi_1, \ldots, \xi_K)$,
%\BEN
%\nabla_{\bm\zeta} F(\bm\xi+\e \tilde{\bm\zeta})   = \e^2 \nabla_{\bm\zeta} \Gamma(\bm\xi)+ \e^{3}\nabla_{\bm\zeta}^2 \Gamma(\bm\xi) \cdot\tilde{\bm\zeta}+ O(\e^4) =\e^{3}\nabla_{\bm\zeta}^2 \Gamma(\bm\xi) \cdot\tilde{\bm\xi}+ O(\e^4).
%\EEN
%On account of the non-degeneracy of $\bm\xi$, Brouwer's fixed point theorem guarantee the existence of $\tilde{\bm\zeta}\in \left(B_{r_0}(0)\right)^K$ such that $\nabla_{\bm\zeta} F(\bm\xi+\e \tilde{\bm\zeta}) =0$.
\end{proof}

%\subsection*{Acknowledgement}

\end{document}